\documentclass[11pt]{article}%
\usepackage{amsmath}
\usepackage{graphicx}%
\usepackage{amsfonts}%
\usepackage{amssymb}
\newcommand{\eqnb}{\begin{equation}}
\newcommand{\eqne}{\end{equation}}

\begin{document}

\title{A Computational Framework for the Mixing Times in the QBD Processes with
Infinitely-Many Levels}
\author{Quan-Lin Li\\School of Economics and Management Sciences \\Yanshan University, Qinhuangdao 066004, P.R. China\\Jing Cao\\College of Science, Yanshan University, Qinhuangdao 066004, P.R. China}
\maketitle

\begin{abstract}
In this paper, we develop some matrix Poisson's equations satisfied by the
mean and variance of the mixing time in an irreducible positive-recurrent
discrete-time Markov chain with infinitely-many levels, and provide a
computational framework for the solution to the matrix Poisson's equations by
means of the UL-type of $RG$-factorization as well as the generalized
inverses. In an important special case: the level-dependent QBD processes, we
provide a detailed computation for the mean and variance of the mixing time.
Based on this, we give new highlight on computation of the mixing time in the
block-structured Markov chains with infinitely-many levels through the
matrix-analytic method.

\vskip                                                     0.3cm


\end{abstract}

\section{Introduction}

In probability theory, the mixing time of a Markov chain is the time until the
Markov chain is close to its steady state distribution. Up to now, the mixing
time has been given many important applications to, for example, perturbation
analysis, Poisson's equations, coupling, spectral gap, random walks on graphs,
and Markov chain Monte Carlo algorithms. During the last two decades
considerable attention has been paid to studying the mixing times in the
Markov chains with finite states. Readers may refer to, such as, Aldous et al
\cite{Ald:1997}, Lov\'{a}sz and Winkler \cite{Lov:1998}, Hunter
\cite{Hun:2006, Hun:2008, Hun:2009, Hun:2013}, Cao et al \cite{Cao:1996}, Cao
and Chen \cite{Cao:1997}, and Li and Liu \cite{Li:2004}. At the same time, the
mixing times in Markov chains and random walks are studied in five excellent
books by Aldous \cite{Ald:1983}, Aldous and Fill \cite{Ald:2002}, Montenegro
and Tetali \cite{Mon:2006}, Cao \cite{Cao:2007}, and Levin et al
\cite{Lev:2009}. It is worthwhile to note that for the mixing times in the
Markov chains with infinite states, the available works are still few in the
literature up to now.

The generalized inverses play an important role in the study of the mixing
times in the Markov chains with finite states. Readers may refer to a book by
Kemeny and Snell \cite{Kem:1960}, and three survey papers by \cite{Mey:1975}
and Hunter \cite{Hun:1982, Hun:1988}. Specifically, Hunter \cite{Hun:1982,
Hun:1988} indicated how to apply the generalized inverses to computing the
mean and variance of the mixing time. On the other hand, the fundamental
matrix plays a key role in the theory of Markov chains, e.g., see Kemeny and
Snell \cite{Kem:1960}, Hunter \cite{Hun:1982}, Neuts \cite{Neu:1989}, Heyman
and O'Leary \cite{Hey:1995}, and da Silva Soares and Latouche \cite{da:2002}.
Also, the kemeny's constant has been an interesting topic for many years,
readers may refer to Hunter \cite{Hun:2011} and Catral et al \cite{Cat:2010}
for more details.

The two types of $RG$-factorizations have been a key method in performance
computation of stochastic models, e.g., see Li \cite{Li:2010} for a
systematical analysis. Vigon \cite{Vig:2013} discussed the LU-factorization
and Wiener-Hopf factorization in Markov chains, and further provided new
highlight on some useful relations among the $RG$-factorizations,
LU-factorization and Wiener-Hopf factorization. From a viewpoint of
applications, the $RG$-factorizations have been applied to dealing with
performance computation in a variety of stochastic models. Important examples
include quasi-stationary distributions by Li and Zhao \cite{Li:2002, Li:2003},
stochastic functionals by Li and Cao \cite{LiC:2004}, tail probabilities by Li
\cite{Li:2013}, repairable systems by Ruiz-Castro et al \cite{Rui:2009},
computer networks by Wang et al \cite{Wan:2007, Wan:2010}, manufacturing
systems by Li et al \cite{Li:2009} and Liu et al \cite{Liu:2012}.

The QBD process is an important example in Markov chains, and provides a
useful mathematical tool for studying stochastic models such as queueing
systems, manufacturing systems, communication networks and transportation
systems. Readers may refer to Chapter 3 of Neuts \cite{Neu:1981}, Bright and
Taylor \cite{Bri:1996}, Ramaswami \cite{Ram:1996}, Latouche and Ramaswami
\cite{Lat:1999}, and Li and Cao \cite{LiC:2004}.

The Poisson's equation frequently occur in the analysis of , such as, Markov
chains, Markov decision processes and queueing systems. Important examples
include the waiting time of a queue by Asmussen and Bladt \cite{Asm:1994},
Glynn \cite{Gly:1994} and Bladt \cite{Bla:1996}; perturbation analysis in
Markov chains by Cao et al \cite{Cao:1996}, Li and Liu \cite{Li:2004}, Cao
\cite{Cao:2007} and Dendievel et al \cite{Den:2013}; Markov decision processes
by Cao \cite{Cao:2007} and Makowski and Shwartz \cite{Mak:2002}. For a
comprehensive analysis of the Poisson's equation, readers may refer to
Nummelin \cite{Num:1991}, Meyn and Tweedie \cite{Mey:1993}, Glynn
\cite{Gly:1996}, Makowski and Shwartz \cite{Mak:2002}.

The main contributions of this paper are twofold. The first one is to develop
some matrix Poisson's equations in the study of mixing times of Markov chains
with either finite states or infinite states. When the Markov chain has finite
states, the generalized inverses are always utilized to study the mean and
variance of the mixing time, e.g., see \cite{Hun:2006, Hun:2008, Hun:2009,
Hun:2013}. However, for a Markov chain with infinite states, the available
works for the solution to the matrix Poisson's equations are few (e.g., see
Nummelin \cite{Num:1991}, Cao and Chen \cite{Cao:1997}, Makowski and Shwartz
\cite{Mak:2002}, Li and Liu \cite{Li:2004} and Dendievel et al \cite{Den:2013}%
), but the mean and variance of the mixing time have not yet been studied in
the literature up to now. This motivates us in this paper to apply the UL-type
of $RG$-factorization as well as the generalized inverses to setting up a
computational framework for solving the matrix Poisson's equations, which can
lead to a detailed analysis for the mixing times in the Markov chains with
infinite states. The second contribution of this paper is to provide a
systematical discussion on the matrix Poisson's equations satisfied by the
mean and variance of the mixing time in the level-dependent QBD process. Our
main results are that the $R$-, $U$- and $G$-measures are used in expressions
both for the solution to the matrix Poisson's equations and for the mean and
variance of the mixing time. Note that some effective algorithms have been
developed for computing the $R$-, $U$- and $G$-measures in the QBD processes
(e.g., see Bright and Taylor \cite{Bri:1995, Bri:1996}), thus this paper
provides effectively numerical computation for the mean and variance of the
mixing time through the matrix-analytic method. On the other hand, the first
passage time described in this paper has a general block-structured
probability meaning, which is different from the first passage times in the
$R$-, $U$- and $G$-measures by means of the taboo probability (e.g., see Neuts
\cite{Neu:1981, Neu:1989} and Li \cite{Li:2010}). Based on this, we develop an
interesting and new research line on which the $R$-, $U$- and $G$-measures are
used to be able to deal with the general first passage time, the mixing time
and the matrix Poission's equations in the Markov chains with infinite states.

The remainder of this paper is organized as follows. In Section 2, for a
general block-structured Markov chain we develop some matrix Poisson's
equations satisfied by the means and variances of the first passage time and
of the mixing time. Furthermore, we apply the UL-type of $RG$-factorization as
well as the generalized inverse to providing a computational procedure for
solving the matrix Poisson's equations. In Section 3, we simply review the
UL-type of $RG$-factorization of the QBD process, and provide explicit
expressions for the two key matrices $\left(  I-R_{U}\right)  ^{-1}$ and
$\left(  I-G_{L}\right)  ^{-1}$ which are given by the $R$- and $G$-measures,
respectively. For the level-dependent QBD process, Section 4 computes the
means of the first passage time and of the mixing time, while Section 5
discusses the variances of the first passage time and of the mixing time. Some
remarks and conclusions are given in the final section.

\section{The Matrix Poisson's Equations}

In this section, for a general block-structured Markov chain with
infinitely-many levels, we develop useful matrix Poisson's equations satisfied
by the means and variances of the first passage time and of the mixing time.
Furthermore, we apply the UL-type of $RG$-factorization as well as the
generalized inverse to providing a computational procedure for solving these
matrix Poisson's equations by means of the $R$- and $G$-measures.

We consider a general discrete-time block-structured Markov chain $\left\{
\left(  X_{n},J_{n}\right)  :n\geq0\right\}  $ whose transition probability
matrix is given by%
\[
P=\left(
\begin{array}
[c]{cccc}%
P_{0,0} & P_{0,1} & P_{0,2} & \cdots\\
P_{1,0} & P_{1,1} & P_{1,2} & \cdots\\
P_{2,0} & P_{2,1} & P_{2,2} & \cdots\\
\vdots & \vdots & \vdots & \ddots
\end{array}
\right)  ,
\]
where $X_{n}$ is the level process and $J_{n}$ is the phase process, the size
of the block $P_{k,k}$ is $m_{k}$ and the sizes of all the other blocks can be
determined accordingly. We assume that the Markov chain $P$ is irreducible and
positive recurrent. In this case, the matrix $P$ is stochastic, that is,
$Pe=e$, $e$ is a column vector of ones with a suitable size. The stationary
probability vector of the Markov chain $P$ is partitioned accordingly into
vectors $\pi=\left(  \pi_{0},\pi_{1},\pi_{2},\ldots\right)  $, where the size
of the vector $\pi_{k}$ is $m_{k}$ for $k\geq0$.

When the Markov chain $\left\{  \left(  X_{n},J_{n}\right)  :n\geq0\right\}  $
is irreducible and positive recurrent, it is clear that%
\[
\underset{n\rightarrow\infty}{\lim}P\left\{  \left(  X_{n},J_{n}\right)
=\left(  k,j\right)  \text{ }|\text{ }\left(  X_{0},J_{0}\right)  =\left(
l,i\right)  \right\}  =\pi_{k,j},
\]
which is independent of the initial state $\left(  X_{0},J_{0}\right)
=\left(  l,i\right)  $. If for some $r\geq0$, $P\left\{  \left(  X_{r}%
,J_{r}\right)  =\left(  k,j\right)  \right\}  =\pi_{k,j}$, then for $n\geq r$%
\[
P\left\{  \left(  X_{n},J_{n}\right)  =\left(  k,j\right)  \right\}
=\pi_{k,j}\text{.}%
\]

\subsection{The mean of the mixing time}

Let $T_{l,i;k,j}$ be the first passage time of the Markov chain $\left\{
\left(  X_{n},J_{n}\right)  :n\geq0\right\}  $ from state $\left(  l,i\right)
$ to state $\left(  k,j\right)  $, that is,%
\[
T_{l,i;k,j}=\min\left\{  n\geq1:\left(  X_{n},J_{n}\right)  =\left(
k,j\right)  \text{ }|\text{ }\left(  X_{0},J_{0}\right)  =\left(  l,i\right)
\right\}  .
\]
Specifically, $T_{k,j;k,j}$ is the first return time of the Markov chain
$\left\{  \left(  X_{n},J_{n}\right)  :n\geq0\right\}  $ to state $\left(
k,j\right)  $.

We write%
\[
M_{l,i;k,j}=E\left[  T_{l,i;k,j}\right]  ,
\]%
\[
M_{l,k}=\left(
\begin{array}
[c]{cccc}%
M_{l,1;k,1} & M_{l,1;k,2} & \cdots & M_{l,1;k,m_{k}}\\
M_{l,2;k,1} & M_{l,2;k,2} & \cdots & M_{l,2;k,m_{k}}\\
\vdots & \vdots &  & \vdots\\
M_{l,m_{l};k,1} & M_{l,m_{l};k,2} & \cdots & M_{l,m_{l};k,m_{k}}%
\end{array}
\right)
\]
and%
\[
\mathbf{M=}\left(
\begin{array}
[c]{cccc}%
M_{0,0} & M_{0,1} & M_{0,2} & \cdots\\
M_{1,0} & M_{1,1} & M_{1,2} & \cdots\\
M_{2,0} & M_{2,1} & M_{2,2} & \cdots\\
\vdots & \vdots & \vdots &
\end{array}
\right)  .
\]
Let $E=ee^{T}$, where $A^{T}$ denotes the transpose of the matrix $A$. We
write%
\[
\text{diag}\left(  \pi_{k}\right)  =\text{diag}\left(  \pi_{k,1},\pi
_{k,2},\ldots,\pi_{k,m_{k}}\right)  ,\text{ \ }k\geq0,
\]
and%
\[
\text{diag}\left(  \pi\right)  =\text{diag}\left(  \text{diag}\left(  \pi
_{0}\right)  ,\text{diag}\left(  \pi_{1}\right)  ,\text{diag}\left(  \pi
_{2}\right)  ,\ldots\right)  .
\]
Then it follows from Theorem 4.4.4 in Kemeny and Snell \cite{Kem:1960} or
(2.1) in Hunter \cite{Hun:2006} that%
\begin{equation}
\left(  I-P\right)  \mathbf{M}=E-P\left[  \text{diag}\left(  \pi\right)
\right]  ^{-1}. \label{FPT-0}%
\end{equation}

When the Markov chain $P$ is at steady-state, the $\left(  \sum_{i=0}%
^{k-1}m_{i}\right.  $ $\left.  +j\right)  $th entry of the column vector
$\mathbf{M}e$ is the mean of stationary first passage time to state $\left(
k,j\right)  $, hence we have%
\[
\tau=\pi\mathbf{M}e
\]
is the mean of stationary first passage time in the Markov chain $\left\{
\left(  X_{n},J_{n}\right)  :n\geq0\right\}  $.

Let $T$ be the mixing time of the Markov chain $\left\{  \left(  X_{n}%
,J_{n}\right)  :n\geq0\right\}  $, and $Y$ a random variable whose probability
distribution is the stationary probability vector $\pi=\left(  \pi_{0},\pi
_{1},\pi_{2},\ldots\right)  $. Then%
\[
T=\min\left\{  n:\left(  X_{n},J_{n}\right)  =Y\right\}  ,
\]
and the Markov chain $\left\{  \left(  X_{n},J_{n}\right)  :n\geq0\right\}  $
reaches stationary or achieves mixing at time $T$.

We define the mean of the mixing time as%
\[
L_{l,i;k,j}=E\left[  T\text{ }|\text{ }\left(  X_{0},J_{0}\right)  =\left(
l,i\right)  ,\left(  X_{T},J_{T}\right)  =\left(  k,j\right)  \right]  .
\]
Let%
\[
L_{l,k}=\left(
\begin{array}
[c]{cccc}%
L_{l,1;k,1} & L_{l,1;k,2} & \cdots & L_{l,1;k,m_{k}}\\
L_{l,2;k,1} & L_{l,2;k,2} & \cdots & L_{l,2;k,m_{k}}\\
\vdots & \vdots &  & \vdots\\
L_{l,m_{l};k,1} & L_{l,m_{l};k,2} & \cdots & L_{l,m_{l};k,m_{k}}%
\end{array}
\right)
\]
and%
\[
\mathbf{L=}\left(
\begin{array}
[c]{cccc}%
L_{0,0} & L_{0,1} & L_{0,2} & \cdots\\
L_{1,0} & L_{1,1} & L_{1,2} & \cdots\\
L_{2,0} & L_{2,1} & L_{2,2} & \cdots\\
\vdots & \vdots & \vdots &
\end{array}
\right)  .
\]
From (2.6) in Hunter \cite{Hun:2006}%
\[
L_{l,i;k,j}=M_{l,i;k,j}\pi_{k,j}.
\]
we obtain%
\begin{equation}
\mathbf{L}=\mathbf{M}\text{diag}\left(  \pi\right)  . \label{MTM-0}%
\end{equation}
Using (\ref{FPT-0}) and (\ref{MTM-0}), we obtain%
\begin{equation}
\left(  I-P\right)  \mathbf{L}=e\pi-P. \label{MTMEqu-0}%
\end{equation}

We write%
\[
\eta_{k,j}=E\left[  T\text{ }|\text{ }\left(  X_{0},J_{0}\right)  =\left(
k,j\right)  \right]  ,
\]%
\[
\eta_{k}=\left(  \eta_{k,1},\eta_{k,2},\ldots,\eta_{k,m_{k}}\right)
\]
and%
\[
\mathbf{\eta}=\left(  \eta_{0},\eta_{1},\eta_{2},\eta_{3},\ldots\right)  .
\]
It follows from Theorem 2.1 in Hunter \cite{Hun:2006} that%
\begin{equation}
\mathbf{\eta}^{T}=\mathbf{M}\pi^{T}, \label{MT-0}%
\end{equation}
Using (\ref{MT-0}) and (\ref{MTM-0}), we obtain%
\[
\mathbf{\eta}^{T}=\mathbf{L}e,
\]
which, together with (\ref{MTMEqu-0}), follows%
\begin{equation}
\left(  I-P\right)  \mathbf{\eta}^{T}=0. \label{MTEqu-0}%
\end{equation}

\subsection{The variance of the mixing time}

Let%
\[
M_{l,i;k,j}^{\left(  2\right)  }=E\left[  T_{l,i;k,j}^{2}\right]  .
\]
We write%
\[
M_{l,k}^{\left(  2\right)  }=\left(
\begin{array}
[c]{cccc}%
M_{l,1;k,1}^{\left(  2\right)  } & M_{l,1;k,2}^{\left(  2\right)  } & \cdots &
M_{l,1;k,m_{k}}^{\left(  2\right)  }\\
M_{l,2;k,1}^{\left(  2\right)  } & M_{l,2;k,2}^{\left(  2\right)  } & \cdots &
M_{l,2;k,m_{k}}^{\left(  2\right)  }\\
\vdots & \vdots &  & \vdots\\
M_{l,m_{l};k,1}^{\left(  2\right)  } & M_{l,m_{l};k,2}^{\left(  2\right)  } &
\cdots & M_{l,m_{l};k,m_{k}}^{\left(  2\right)  }%
\end{array}
\right)
\]
and%
\[
\mathbf{M}^{\left(  2\right)  }\mathbf{=}\left(
\begin{array}
[c]{cccc}%
M_{0,0}^{\left(  2\right)  } & M_{0,1}^{\left(  2\right)  } & M_{0,2}^{\left(
2\right)  } & \cdots\\
M_{1,0}^{\left(  2\right)  } & M_{1,1}^{\left(  2\right)  } & M_{1,2}^{\left(
2\right)  } & \cdots\\
M_{2,0}^{\left(  2\right)  } & M_{2,1}^{\left(  2\right)  } & M_{2,2}^{\left(
2\right)  } & \cdots\\
\vdots & \vdots & \vdots &
\end{array}
\right)  .
\]

Using (2.7) in Hunter \cite{Hun:2008}, we obtain%
\begin{equation}
\left(  I-P\right)  \mathbf{M}^{\left(  2\right)  }=E+P\left\{  2\mathbf{M}%
--\left[  \text{diag}\left(  \pi\right)  \right]  ^{-1}\left[  I+2\left(
e\pi\mathbf{M}\right)  _{d}\right]  \right\}  ,\label{M2Equ-0}%
\end{equation}
where $\left(  e\pi\mathbf{M}\right)  _{d}$ is a diagonal matrix whose
diagonal entries are given by the diagonal entries of the matrix
$e\pi\mathbf{M}$.

Let%
\[
L_{l,i;k,j}^{\left(  2\right)  }=E\left[  T^{2}\text{ }|\text{ }\left(
X_{0},J_{0}\right)  =\left(  l,i\right)  ,\left(  X_{T},J_{T}\right)  =\left(
k,j\right)  \right]  .
\]
Then%
\begin{equation}
L_{l,i;k,j}^{\left(  2\right)  }=M_{l,i;k,j}^{\left(  2\right)  }\pi
_{k,j}.\label{L2Pi}%
\end{equation}
We write%
\[
L_{l,k}^{\left(  2\right)  }=\left(
\begin{array}
[c]{cccc}%
L_{l,1;k,1}^{\left(  2\right)  } & L_{l,1;k,2}^{\left(  2\right)  } & \cdots &
L_{l,1;k,m_{k}}^{\left(  2\right)  }\\
L_{l,2;k,1}^{\left(  2\right)  } & L_{l,2;k,2}^{\left(  2\right)  } & \cdots &
L_{l,2;k,m_{k}}^{\left(  2\right)  }\\
\vdots & \vdots &  & \vdots\\
L_{l,m_{l};k,1}^{\left(  2\right)  } & L_{l,m_{l};k,2}^{\left(  2\right)  } &
\cdots & L_{l,m_{l};k,m_{k}}^{\left(  2\right)  }%
\end{array}
\right)
\]
and%
\[
\mathbf{L}^{\left(  2\right)  }\mathbf{=}\left(
\begin{array}
[c]{cccc}%
L_{0,0}^{\left(  2\right)  } & L_{0,1}^{\left(  2\right)  } & L_{0,2}^{\left(
2\right)  } & \cdots\\
L_{1,0}^{\left(  2\right)  } & L_{1,1}^{\left(  2\right)  } & L_{1,2}^{\left(
2\right)  } & \cdots\\
L_{2,0}^{\left(  2\right)  } & L_{2,1}^{\left(  2\right)  } & L_{2,2}^{\left(
2\right)  } & \cdots\\
\vdots & \vdots & \vdots &
\end{array}
\right)  .
\]
It is easy to see from (\ref{L2Pi}) that%
\begin{equation}
\mathbf{L}^{\left(  2\right)  }=\mathbf{M}^{\left(  2\right)  }\text{diag}%
\left(  \pi\right)  .\label{L2}%
\end{equation}
It follows from (\ref{M2Equ-0}) that%
\begin{equation}
\left(  I-P\right)  \mathbf{L}^{\left(  2\right)  }=e\pi+P\left\{
2\mathbf{M}-\left[  \text{diag}\left(  \pi\right)  \right]  ^{-1}\left[
I+2\left(  e\pi\mathbf{M}\right)  _{d}\right]  \right\}  \text{diag}\left(
\pi\right)  .\label{L2Equ-0}%
\end{equation}

Set%
\[
\eta_{k,j}^{\left(  2\right)  }=E\left[  T^{2}\text{ }|\text{ }\left(
X_{0},J_{0}\right)  =\left(  k,j\right)  \right]  ,
\]%
\[
\eta_{k}^{\left(  2\right)  }=\left(  \eta_{k,1}^{\left(  2\right)  }%
,\eta_{k,2}^{\left(  2\right)  },\ldots,\eta_{k,m_{k}}^{\left(  2\right)
}\right)
\]
and%
\[
\mathbf{\eta}^{\left(  2\right)  }=\left(  \eta_{0}^{\left(  2\right)  }%
,\eta_{1}^{\left(  2\right)  },\eta_{2}^{\left(  2\right)  },\eta_{3}^{\left(
2\right)  },\ldots\right)  .
\]
It follows from Theorem 1.1 in Hunter \cite{Hun:2008} that%
\begin{equation}
\left(  \mathbf{\eta}^{\left(  2\right)  }\right)  ^{T}=\mathbf{M}^{\left(
2\right)  }\pi^{T},\label{eta2-0}%
\end{equation}
Based on (\ref{M2Equ-0}) and (\ref{eta2-0}), we obtain%
\begin{equation}
\left(  I-P\right)  \left(  \mathbf{\eta}^{\left(  2\right)  }\right)
^{T}=e+P\left\{  2\mathbf{M}-\left[  \text{diag}\left(  \pi\right)  \right]
^{-1}\left[  I+2\left(  e\pi\mathbf{M}\right)  _{d}\right]  \right\}  \pi
^{T}.\label{eta2Equ-0}%
\end{equation}

Let%
\[
v_{k,j}^{\left(  2\right)  }=Var\left[  T\text{ }|\text{ }\left(  X_{0}%
,J_{0}\right)  =\left(  k,j\right)  \right]  ,
\]%
\[
v_{k}^{\left(  2\right)  }=\left(  v_{k,1}^{\left(  2\right)  },v_{k,2}%
^{\left(  2\right)  },\ldots,v_{k,m_{k}}^{\left(  2\right)  }\right)
\]
and%
\[
\mathbf{V}^{\left(  2\right)  }=\left(  v_{0}^{\left(  2\right)  }%
,v_{1}^{\left(  2\right)  },v_{2}^{\left(  2\right)  },v_{3}^{\left(
2\right)  },\ldots\right)  .
\]
Then%
\begin{equation}
\left(  \mathbf{V}^{\left(  2\right)  }\right)  ^{T}=\left(  \mathbf{\eta
}^{\left(  2\right)  }\right)  ^{T}-\mathbf{\eta}^{T}\circ\mathbf{\eta}^{T}
\label{VarV}%
\end{equation}
where $A\circ B$ denotes the Hadamard Product of the two matrices $A$ and $B$,
that is, $A\circ B=\left(  a_{i,j}b_{i,j}\right)  $ if $A=\left(
a_{i,j}\right)  $ and $B=\left(  b_{i,j}\right)  $. Let $v^{\left(  2\right)
}=Var\left(  T\right)  $. Then%
\begin{equation}
v^{\left(  2\right)  }=\pi\left(  \mathbf{V}^{\left(  2\right)  }\right)
^{T}=\pi\left(  \mathbf{\eta}^{\left(  2\right)  }\right)  ^{T}-\pi\left(
\mathbf{\eta}^{T}\circ\mathbf{\eta}^{T}\right)  , \label{Var-0}%
\end{equation}
which is the variance of the mixing time when the Markov chain $P$ is at
steady state.

It is worthwhile to note that the above matrix Poisson's equations are a
direct generalization of those in Hunter \cite{Hun:2008} both from the finite
states to the infinite states, and from the scale entries to the block
entries. Therefore, the matrix Poisson's equations developed here are more
general than those in Hunter \cite{Hun:2008}, and they are useful and
interesting in the study of stochastic models through the matrix-analytic method.

\subsection{The Poisson's equation}

In the UL-type of $RG$-factorization%
\[
I-P=\left(  I-R_{U}\right)  \left(  I-\Psi_{D}\right)  \left(  I-G_{L}\right)
,
\]
where $\Psi_{D}=$ diag$\left(  U_{0},U_{1},U_{2},\ldots\right)  $ and $U_{0}$
is the transition probability matrix of the censoring Markov chain to level 0.
If the Markov chain $P$ is irreducible and recurrent, then the censoring
Markov chain $U_{0}$ is irreducible and positive recurrent, and rank$\left(
U_{0}\right)  =m_{0}-1$. Let $v_{0}$ be the stationary probability vector of
the censoring Markov chain $U_{0}$. Then $v_{0}\left(  I-U_{0}\right)  =0$.

Now, we deal with the Poisson's equation:%
\begin{equation}
\left(  I-U_{0}\right)  x=g, \label{Poisson-1}%
\end{equation}
where $g$ is a given column vector of size $m_{0}$. Hence it follows from
(\ref{Poisson-1}) that $v_{0}g=0$. This shows that the given vector $g$ must
satisfy the condition $v_{0}g=0$ if there exists one solution to the Poisson's
equation (\ref{Poisson-1}).

If there exists a matrix $V$ such that $\left(  I-U_{0}\right)  V\left(
I-U_{0}\right)  =\left(  I-U_{0}\right)  $, then the matrix $V$ is called a
generalized inverse of the matrix $I-U_{0}$.

From Theorem 3.3 in Hunter \cite{Hun:1982}, we know that the matrix
$I-U_{0}+\mathbf{tu}$ is invertible and $\left[  I-U_{0}+\mathbf{tu}\right]
^{-1}$ is a generalized inverse of the matrix $I-U_{0}$, where $\mathbf{t}$
and $\mathbf{u}$ are two arbitrary vectors such that $v_{0}\mathbf{t}\neq0$
and $\mathbf{u}e\neq0$. It is clear that $\mathbf{t}$ and $\mathbf{u}$ are the
column and row vectors, respectively. Furthermore, the matrix%
\begin{equation}
V=\left[  I-U_{0}+\mathbf{tu}\right]  ^{-1}+e\mathbf{f}+\mathbf{h}v_{0}
\label{Poisson-2}%
\end{equation}
is a generalized inverse of the matrix $I-U_{0}$, where $\mathbf{f}$ and
$\mathbf{h}$ are two arbitrary vectors. Clearly, $\mathbf{f}$ and $\mathbf{h}$
are the row and column vectors, respectively. It is worthwhile to note that
any generalized inverse of the matrix $I-U_{0}$ can be expressed by
(\ref{Poisson-2}) with the four vectors: $\mathbf{t}$, $\mathbf{u}$,
$\mathbf{f}$ and $\mathbf{h}$. Specifically, when $\mathbf{t}=e$,
$\mathbf{u}=v_{0}$, $\mathbf{f}=0$ and $\mathbf{h}=0$, the matrix%
\begin{equation}
Z=\left[  I-U_{0}+ev_{0}\right]  ^{-1} \label{Poisson-3}%
\end{equation}
is called the Kemeny and Snell's fundamental matrix.

It is easy to check that%
\[
\mathbf{u}\left[  I-U_{0}+\mathbf{tu}\right]  ^{-1}=\frac{v_{0}}%
{v_{0}\mathbf{t}}%
\]
and%
\[
\left[  I-U_{0}+\mathbf{tu}\right]  ^{-1}\mathbf{t}=\frac{e}{\mathbf{u}e}.
\]
Specifically, we have%
\[
v_{0}\left[  I-U_{0}+ev_{0}\right]  ^{-1}=v_{0}%
\]
and%
\[
\left[  I-U_{0}+ev_{0}\right]  ^{-1}e=e.
\]

Let $V$ be any generalized inverse of the matrix $I-U_{0}$, given by
(\ref{Poisson-2}). Then the Poisson's equation $\left(  I-U_{0}\right)  x=g$
has a solution if and only if%
\[
\left(  I-U_{0}\right)  Vg=g.
\]
In this case, we have
\begin{equation}
x=Vg+\left[  I-V\left(  I-U_{0}\right)  \right]  \Theta, \label{Poisson-4}%
\end{equation}
where $\Theta$ is an arbitrary column vector. Specifically, we have an
important solution as follows:%
\begin{equation}
x=\left[  I-U_{0}+ev_{0}\right]  ^{-1}g+ce, \label{Poisson-5}%
\end{equation}
where $c$ is an arbitrary constant.

Now, we give some useful observation on the solution (\ref{Poisson-5}) as
follows:%
\begin{align*}
\left(  I-U_{0}\right)  x  &  =\left(  I-U_{0}\right)  \left[  I-U_{0}%
+ev_{0}\right]  ^{-1}g+\left(  I-U_{0}\right)  ce\\
&  =g-ev_{0}g=g,
\end{align*}
since $\left(  I-U_{0}\right)  e=0$ and the basic condition $v_{0}g=0$.

Because of the fact that rank$\left(  U_{0}\right)  =m_{0}-1$, it seems to
have a better understanding for the only one free parameter $c$ in the
solution (\ref{Poisson-5}). On the contrary, the solution (\ref{Poisson-4})
have more free parameters through the five vectors $\mathbf{t}$, $\mathbf{u}$,
$\mathbf{f}$, $\mathbf{h}$ and $\Theta$.

In the rest of this subsection, we consider a matrix Poisson's equation%
\begin{equation}
\left(  I-U_{0}\right)  \mathbf{X}=\mathbf{G}, \label{Poisson-6}%
\end{equation}
where $\mathbf{G}$ is a given matrix of size $m_{0}$. To solve the matrix
Poisson's equation, we write%
\[
\mathbf{X}=\left(  x_{1},x_{2},\ldots,x_{m_{0}}\right)
\]
and%
\[
\mathbf{G}=\left(  g_{1},g_{2},\ldots,g_{m_{0}}\right)  .
\]
Then the matrix Poisson's equation is written as the $m_{0}$ Poisson's
equations as follows:%
\[
\left(  I-U_{0}\right)  x_{k}=g_{k},\text{ \ }1\leq k\leq m_{0}.
\]
It is easy to see from (\ref{Poisson-5}) and (\ref{Poisson-6}) that%
\begin{equation}
\mathbf{X}=\left[  I-U_{0}+ev_{0}\right]  ^{-1}\mathbf{G}+e\mathbf{c},
\label{Poisson-7}%
\end{equation}
where $\mathbf{c}=\left(  c_{1},c_{2},\ldots,c_{m_{0}}\right)  $ are an
arbitrary row vector. It is seen that the solution (\ref{Poisson-7}) contains
the most basic $m_{0}$ free parameters in the vector $\mathbf{c}=\left(
c_{1},c_{2},\ldots,c_{m_{0}}\right)  $. Note that the solution
(\ref{Poisson-7}) will be useful in the remainder of this paper.

\subsection{A computational procedure}

From the above matrix Poisson's equations (e.g., (\ref{FPT-0}) and
(\ref{M2Equ-0})), it is seen that we need to solve a general matrix Poisson's
equation as follows:%
\begin{equation}
\left(  I-P\right)  \mathbf{A}=\mathbf{B}, \label{MatrixE-1}%
\end{equation}
where $\mathbf{B}$ is a given matrix with $\pi\mathbf{B}=0$. From Chapter 2 in
Li \cite{Li:2010}, the UL-type of $RG$-factorization for the general Markov
chain $P$ is given by%
\[
I-P=\left(  I-R_{U}\right)  \left(  I-\Psi_{D}\right)  \left(  I-G_{L}\right)
.
\]
This gives%
\[
\left(  I-R_{U}\right)  \left(  I-\Psi_{D}\right)  \left(  I-G_{L}\right)
\mathbf{A}=\mathbf{B},
\]
which follows%
\[
\left(  I-\Psi_{D}\right)  \left(  I-G_{L}\right)  \mathbf{A}=\left(
I-R_{U}\right)  ^{-1}\mathbf{B},
\]

Let%
\begin{equation}
\mathbf{X}=\left(  I-G_{L}\right)  \mathbf{A}. \label{MatrixE-2}%
\end{equation}
Then%
\begin{equation}
\left(  I-\Psi_{D}\right)  \mathbf{X=}\left(  I-R_{U}\right)  ^{-1}\mathbf{B}.
\label{MatrixE-3}%
\end{equation}
We write%
\[
I-\Psi_{D}=\text{diag}\left(  I-U_{0},I-\Phi_{D}\right)  ,
\]%
\[
\Phi_{D}=\text{diag}\left(  U_{1},U_{2},U_{3},\ldots\right)  ;
\]%
\[
\mathbf{X}=\left(
\begin{array}
[c]{c}%
\mathbf{X}_{0}\\
\mathbf{X}_{1}%
\end{array}
\right)  ,
\]
where $\mathbf{X}_{0}$ is a matrix with the first $m_{0}$ row vectors of the
matrix $\mathbf{X}$;%
\begin{equation}
\left(  I-R_{U}\right)  ^{-1}\mathbf{B}=\left(
\begin{array}
[c]{c}%
\mathbf{C}_{0}\\
\mathbf{C}_{1}%
\end{array}
\right)  , \label{MatrixE-3-1}%
\end{equation}
where $\mathbf{C}_{0}$ is a matrix with the first $m_{0}$ row vectors of the
matrix $\left(  I-R_{U}\right)  ^{-1}\mathbf{B}$. It follows from
(\ref{MatrixE-2}) that%
\begin{equation}
\left(  I-U_{0}\right)  \mathbf{X}_{0}=\mathbf{C}_{0} \label{MatrixE-4}%
\end{equation}
and%
\begin{equation}
\left(  I-\Phi_{D}\right)  \mathbf{X}_{1}=\mathbf{C}_{1}. \label{MatrixE-5}%
\end{equation}
Note that the matrix $I-U_{k}$ is invertible for $k\geq1$, the matrix
$I-\Phi_{D}$ is invertible. Thus it follows from (\ref{MatrixE-5}) that%
\begin{equation}
\mathbf{X}_{1}=\left(  I-\Phi_{D}\right)  ^{-1}\mathbf{C}_{1}.
\label{MatrixE-6}%
\end{equation}

Note that the matrix $I-U_{0}$ is singular, thus it follows from
(\ref{MatrixE-4}) and (\ref{Poisson-7}) that%
\begin{equation}
\mathbf{X}_{0}=Z\mathbf{C}_{0}+ec_{0}, \label{MatrixE-7}%
\end{equation}
where $c_{0}$ is an arbitrary row vector.

Based on (\ref{MatrixE-6}) and (\ref{MatrixE-7}), we obtain%
\begin{equation}
\mathbf{X}=\left(
\begin{array}
[c]{c}%
\mathbf{X}_{0}\\
\mathbf{X}_{1}%
\end{array}
\right)  =\left(
\begin{array}
[c]{c}%
Z\mathbf{C}_{0}+ec_{0}\\
\left(  I-\Phi_{D}\right)  ^{-1}\mathbf{C}_{1}%
\end{array}
\right)  . \label{MatrixE-8}%
\end{equation}
It follows from (\ref{MatrixE-2}) that%
\[
\left(  I-G_{L}\right)  \mathbf{A=}\left(
\begin{array}
[c]{c}%
Z\mathbf{C}_{0}+ec_{0}\\
\left(  I-\Phi_{D}\right)  ^{-1}\mathbf{C}_{1}%
\end{array}
\right)  .
\]
Thus we obtain%
\begin{equation}
\mathbf{A=}\left(  I-G_{L}\right)  ^{-1}\left(
\begin{array}
[c]{c}%
Z\mathbf{C}_{0}+ec_{0}\\
\left(  I-\Phi_{D}\right)  ^{-1}\mathbf{C}_{1}%
\end{array}
\right)  . \label{MatrixE-9}%
\end{equation}

From (\ref{MatrixE-9}) and (\ref{MatrixE-3-1}), it is seen that a key for the
solution to the matrix Poisson's equation (\ref{MatrixE-1}) is that the two
matrices $\left(  I-R_{U}\right)  ^{-1}$ and $\left(  I-G_{L}\right)  ^{-1}$
can have the explicit expressions by means of the $R$- and $G$-measures.
However, Chapter 2 and Appendix B in Li \cite{Li:2010} indicated that $\left(
I-R_{U}\right)  ^{-1}$ and $\left(  I-G_{L}\right)  ^{-1}$ can explicitly be
expressed only in two special cases: The QBD processes, and Markov chains of
GI/G/1 type.

Note that the first passage time in the matrix $\mathbf{M}$ is different from
the first passage times in the $R$- and $G$-measures, because the $R$- and
$G$-measures are defined by the taboo probability, e.g., see Chapter 2 in Li
\cite{Li:2010} and Neuts \cite{Neu:1981, Neu:1989}. Therefore, this paper
provides new highlight on the block-structured Markov chains, including the
QBD processes, and the Markov chains of GI/M/1 type and of M/G/1 type.

In the remainder of this paper, we will consider an important case: the QBD
processes, and derive the means and variances of the first passage time and of
the mixing time.

\section{The QBD Processes}

In this section, we consider an irreducible discrete-time level-dependent QBD
process with infinitely-many levels, and simply review the UL-type of
$RG$-factorization of the QBD process. Specifically, we provide explicit
expressions for the two key matrices $\left(  I-R_{U}\right)  ^{-1}$ and
$\left(  I-G_{L}\right)  ^{-1}$ through the $R$- and $G$-measures, respectively.

We consider an irreducible discrete-time level-dependent QBD process $\left\{
\left(  X_{n},J_{n}\right)  ,n\geq0\right\}  $ whose transition probability
matrix is given by
\begin{equation}
P=\left(
\begin{array}
[c]{ccccc}%
A_{1}^{(0)} & A_{0}^{(0)} &  &  & \\
A_{2}^{(1)} & A_{1}^{(1)} & A_{0}^{(1)} &  & \\
& A_{2}^{(2)} & A_{1}^{(2)} & A_{0}^{(2)} & \\
&  & \ddots & \ddots & \ddots
\end{array}
\right)  , \label{TPM}%
\end{equation}
where the size of the matrix $A_{1}^{\left(  k\right)  }$ is $m_{k}$ for
$k\geq0$, $A_{i}^{(0)}\gvertneqq0$ for $i=0,1$, $A_{j}^{(k)}\gvertneqq0$ for
$j=0,1,2$ and $k\geq1$, $A_{0}^{(0)}e+A_{1}^{(0)}e=e$ and $A_{0}^{(k)}%
e+A_{1}^{(k)}e+A_{2}^{(k)}e=e$ for $k\geq1$.

Let the matrix sequences $\left\{  R_{l}:\text{\ }l\geq0\right\}  $ and
$\left\{  G_{k}:k\geq1\right\}  $ be the minimal nonnegative solutions to the
systems of matrix Poisson's equations%
\begin{equation}
A_{0}^{(l)}+R_{l}A_{1}^{(l+1)}+R_{l}R_{l+1}A_{2}^{(l+2)}=R_{l},\text{
\ \ }l\geq0, \label{Reqn}%
\end{equation}
and%
\begin{equation}
A_{0}^{(k)}G_{k+1}G_{k}+A_{1}^{(k)}G_{k}+A_{2}^{(k)}=G_{k},\text{ \ \ }k\geq1,
\label{Geqn}%
\end{equation}
respectively. Then the $U$-measure $\left\{  U_{l}:\text{\ }l\geq0\right\}  $
is given by%
\begin{equation}
U_{l}=A_{1}^{(l)}+R_{l}A_{2}^{(l+1)}=A_{1}^{(l)}+A_{0}^{(l)}G_{l+1},\text{
\ \ }l\geq0, \label{Uequ}%
\end{equation}
For $k\geq1$, the matrix $I-U_{k}$ is invertible; while the matrix $I-U_{0}$
is singular, and the censoring Markov chain $U_{0}$ is irreducible and
positive\ only if the QBD process $P$ is irreducible and recurrent.

For the QBD process with infinitely-many levels given in (\ref{TPM}), the
UL-type of $RG$-factorization is given by%
\begin{equation}
I-P=\left(  I-R_{U}\right)  \left(  I-\Psi_{D}\right)  \left(  I-G_{L}\right)
, \label{INFULRG}%
\end{equation}
where%
\[
R_{U}=\left(
\begin{array}
[c]{ccccc}%
0 & R_{0} &  &  & \\
& 0 & R_{1} &  & \\
&  & 0 & R_{2} & \\
&  &  & \ddots & \ddots
\end{array}
\right)  ,
\]%
\[
\Psi_{D}=\text{diag}\left(  U_{0},U_{1},U_{2},U_{3},\ldots\right)
\]
and%
\[
\text{ \ }G_{L}=\left(
\begin{array}
[c]{ccccc}%
0 &  &  &  & \\
G_{1} & 0 &  &  & \\
& G_{2} & 0 &  & \\
&  & G_{3} & 0 & \\
&  &  & \ddots & \ddots
\end{array}
\right)  .
\]
Note that the $RG$-factorization can be given only if the QBD process is irreducible.

Let%
\begin{equation}
X_{k}^{\left(  l\right)  }=R_{l}R_{l+1}R_{l+2}\cdots R_{l+k-1},\text{ \ }%
k\geq1,l\geq0, \label{Xkl}%
\end{equation}
and%
\begin{equation}
Y_{k}^{\left(  l\right)  }=G_{l}G_{l-1}G_{l-2}\cdots G_{l-k+1},\text{ \ }l\geq
k\geq1. \label{Ykl}%
\end{equation}
Then%
\begin{equation}
\left(  I-R_{U}\right)  ^{-1}=\left(
\begin{array}
[c]{ccccc}%
I & X_{1}^{\left(  0\right)  } & X_{2}^{\left(  0\right)  } & X_{3}^{\left(
0\right)  } & \cdots\\
& I & X_{1}^{\left(  1\right)  } & X_{2}^{\left(  1\right)  } & \cdots\\
&  & I & X_{1}^{\left(  2\right)  } & \cdots\\
&  &  & I & \cdots\\
&  &  &  & \ddots
\end{array}
\right)  \label{RInverse}%
\end{equation}
and%
\begin{equation}
\left(  I-G_{L}\right)  ^{-1}=\left(
\begin{array}
[c]{ccccc}%
I &  &  &  & \\
Y_{1}^{\left(  1\right)  } & I &  &  & \\
Y_{2}^{\left(  2\right)  } & Y_{1}^{\left(  2\right)  } & I &  & \\
Y_{3}^{\left(  3\right)  } & Y_{2}^{\left(  3\right)  } & Y_{1}^{\left(
3\right)  } & I & \\
\vdots & \vdots & \vdots & \vdots & \ddots
\end{array}
\right)  . \label{GInverse}%
\end{equation}

If the QBD process is level-independent, then%
\[
R_{k}=R,\text{ \ }k\geq1,
\]
and%
\[
G_{l}=G,\text{ \ }l\geq2.
\]
In this case, we obtain%
\begin{equation}
\left(  I-R_{U}\right)  ^{-1}=\left(
\begin{array}
[c]{ccccc}%
I & R_{0} & R_{0}R & R_{0}R^{2} & \cdots\\
& I & R & R^{2} & \cdots\\
&  & I & R & \cdots\\
&  &  & I & \cdots\\
&  &  &  & \ddots
\end{array}
\right)  \label{RInverse-1}%
\end{equation}
and%
\begin{equation}
\left(  I-G_{L}\right)  ^{-1}=\left(
\begin{array}
[c]{ccccc}%
I &  &  &  & \\
G_{1} & I &  &  & \\
GG_{1} & G & I &  & \\
G^{2}G_{1} & G^{2} & G & I & \\
\vdots & \vdots & \vdots & \vdots & \ddots
\end{array}
\right)  . \label{GInverse-1}%
\end{equation}

It is worthwhile to note that we can obtain the solution to the matrix
Poisson's equation (\ref{MatrixE-1}) once the two key matrices $\left(
I-R_{U}\right)  ^{-1}$ and $\left(  I-G_{L}\right)  ^{-1}$ are expressed
explicitly by means of the $R$-measure $\left\{  R_{l}:\text{\ }%
l\geq0\right\}  $ and the $G$-measure $\left\{  G_{k}:k\geq1\right\}  $,
respectively. It is seen from (\ref{RInverse}) to (\ref{GInverse-1}) that the
first passage time as well as the mixing time in the QBD processes can be
given a detailed analysis.

In the remainder of this paper, we assume that the QBD process is irreducible
and positive recurrent. It follows from Li and Cao \cite{Li:2004} or Li
\cite{Li:2010} that
\begin{equation}
\pi_{0}=\varphi v_{0} \label{Vector1}%
\end{equation}
and
\begin{equation}
\pi_{k}=\varphi v_{0}R_{0}R_{1}\cdots R_{k-1},\text{ \ }k\geq1,
\label{Vector2}%
\end{equation}
where $v_{0}$ is the stationary probability vector of the censored Markov
chain $U_{0}=A_{1}^{\left(  0\right)  }+R_{0}A_{2}^{\left(  1\right)  }$ to
level $0$, and the normalized constant $\varphi$ is given by%
\[
\varphi=\frac{1}{1+\sum\limits_{k=0}^{\infty}v_{0}R_{0}R_{1}\cdots R_{k}e}.
\]

\section{The Mean of the Mixing Time}

In this section, for the QBD process we apply the UL-type of $RG$%
-factorization as well as the generalized inverse to computing the means of
the first passage time and of the mixing time. Note that they can be expressed
by means of the $R$- and $G$-measures through the computational procedure
given in Subsection 2.4.

\subsection{The mean of the first passage time}

To compute the matrix $\mathbf{M}$, we need to solve the equation
(\ref{FPT-0}). Using the UL-type of $RG$-factorization, we obtain%
\[
\left(  I-R_{U}\right)  \left(  I-\Psi_{D}\right)  \left(  I-G_{L}\right)
\mathbf{M}=E-P\left[  \text{diag}\left(  \pi\right)  \right]  ^{-1},
\]
which follows%
\begin{equation}
\left(  I-\Psi_{D}\right)  \left(  I-G_{L}\right)  \mathbf{M}=\left(
I-R_{U}\right)  ^{-1}\left\{  E-P\left[  \text{diag}\left(  \pi\right)
\right]  ^{-1}\right\}  .\label{FPT-1}%
\end{equation}
Let%
\begin{equation}
\mathbf{X}=\left(  I-G_{L}\right)  \mathbf{M}.\label{FPT-2}%
\end{equation}
Then%
\begin{equation}
\left(  I-\Psi_{D}\right)  \mathbf{X=}\left(  I-R_{U}\right)  ^{-1}\left\{
E-P\left[  \text{diag}\left(  \pi\right)  \right]  ^{-1}\right\}
.\label{FPT-3}%
\end{equation}
We write%
\[
\mathbf{X}=\left(
\begin{array}
[c]{c}%
\mathbf{X}_{0}\\
\mathbf{X}_{1}%
\end{array}
\right)  ,
\]
where $\mathbf{X}_{0}$ is a matrix with the first $m_{0}$ row vectors of the
matrix $\mathbf{X}$;%
\[
\left(  I-R_{U}\right)  ^{-1}\left\{  E-P\left[  \text{diag}\left(
\pi\right)  \right]  ^{-1}\right\}  =\left(
\begin{array}
[c]{c}%
\mathbf{F}_{0}\\
\mathbf{F}_{1}%
\end{array}
\right)  ,
\]
where%
\[
\mathbf{F}_{0}=\left(  F_{0,0},F_{0,1},F_{0,2},\ldots\right)  ,
\]%
\[
\mathbf{F}_{1}=\left(
\begin{array}
[c]{cccc}%
F_{1,0} & F_{1,1} & F_{1,2} & \cdots\\
F_{2,0} & F_{2,1} & F_{2,2} & \cdots\\
F_{3,0} & F_{3,1} & F_{3,2} & \cdots\\
\vdots & \vdots & \vdots &
\end{array}
\right)  ,
\]%
\[
F_{i,j}=\left\{
\begin{array}
[c]{ll}%
E_{i,j}+\sum\limits_{k=1}^{\infty}X_{k}^{\left(  i\right)  }E_{k+i,j}%
-X_{j-i}^{\left(  i\right)  }\left[  \text{diag}\left(  \pi_{j}\right)
\right]  ^{-1}, & j\geq i+1,i\geq0,\\
E_{i,j}+\sum\limits_{k=1}^{\infty}X_{k}^{\left(  i\right)  }E_{k+i,j}-\left[
A_{1}^{\left(  i\right)  }+X_{1}^{\left(  i\right)  }A_{2}^{\left(
i+1\right)  }\right]  \left[  \text{diag}\left(  \pi_{j}\right)  \right]
^{-1}, & j=i,i\geq0,\\
E_{i,j}+\sum\limits_{k=1}^{\infty}X_{k}^{\left(  i\right)  }E_{k+i,j}%
-A_{2}^{\left(  i\right)  }\left[  \text{diag}\left(  \pi_{j}\right)  \right]
^{-1}, & j=i-1,i\geq1,\\
E_{i,j}+\sum\limits_{k=1}^{\infty}X_{k}^{\left(  i\right)  }E_{k+i,j}, & j\leq
i-2,i\geq2.
\end{array}
\right.
\]

It follows from (\ref{FPT-3}) that%
\begin{equation}
\left(  I-U_{0}\right)  \mathbf{X}_{0}=\mathbf{F}_{0} \label{FPT-4}%
\end{equation}
and%
\begin{equation}
\left(  I-\Phi_{D}\right)  \mathbf{X}_{1}=\mathbf{F}_{1}. \label{FPT-5}%
\end{equation}
It follows from (\ref{FPT-5}) that%
\begin{equation}
\mathbf{X}_{1}=\left(  I-\Phi_{D}\right)  ^{-1}\mathbf{F}_{1}, \label{FPT-6}%
\end{equation}
and from (\ref{FPT-4}) that%
\begin{equation}
\mathbf{X}_{0}=Z\mathbf{F}_{0}+ec_{0}, \label{FPT-7}%
\end{equation}
where $c_{0}=\left(  c_{0,0},c_{0,1},c_{0,2},\ldots\right)  $ is an arbitrary
row vector.

Based on (\ref{FPT-6}) and (\ref{FPT-7}), we obtain%
\begin{equation}
\mathbf{X}=\left(
\begin{array}
[c]{c}%
\mathbf{X}_{0}\\
\mathbf{X}_{1}%
\end{array}
\right)  =\left(
\begin{array}
[c]{c}%
Z\mathbf{F}_{0}+ec_{0}\\
\left(  I-\Phi_{D}\right)  ^{-1}\mathbf{F}_{1}%
\end{array}
\right)  .\label{FPT-8}%
\end{equation}
It follows from (\ref{FPT-2}) that%
\[
\left(  I-G_{L}\right)  \mathbf{M=}\left(
\begin{array}
[c]{c}%
Z\mathbf{F}_{0}+ec_{0}\\
\left(  I-\Phi_{D}\right)  ^{-1}\mathbf{F}_{1}%
\end{array}
\right)  .
\]
Thus we obtain%
\[
\mathbf{M}=\left(  I-G_{L}\right)  ^{-1}\left(
\begin{array}
[c]{c}%
Z\mathbf{F}_{0}+ec_{0}\\
\left(  I-\Phi_{D}\right)  ^{-1}\mathbf{F}_{1}%
\end{array}
\right)  .
\]
This gives%
\[
M_{i,j}=\left\{
\begin{array}
[c]{ll}%
Z\mathbf{F}_{0}+ec_{0,j}, & i=0,j\geq0,\\
Y_{i}^{\left(  i\right)  }\left(  Z\mathbf{F}_{0}+ec_{0,j}\right)  +\left(
I-U_{i}\right)  ^{-1}F_{i,j}+\sum\limits_{k=1}^{i-1}Y_{i-k}^{\left(  i\right)
}\left(  I-U_{k}\right)  ^{-1}F_{k,j}, & i\geq1,j\geq0.
\end{array}
\right.
\]

\subsection{The mean of the mixing time}

From (\ref{MTM-0}), we have $\mathbf{L}=\mathbf{M}$diag$\left(  \pi\right)  $.
Once the matrix $\mathbf{M}$ is given in Subsection 4.1, it is clear that the
matrix $\mathbf{L}$ is obtained by $\mathbf{M}$diag$\left(  \pi\right)  $.

On the other hand, the matrix $\mathbf{L}$ can be solved by the matrix
equation (\ref{MTMEqu-0}) of itself. It may be necessary to simply provide a
simple outline for the solution to Equation (\ref{MTMEqu-0}).

To compute the matrix $\mathbf{L}$ from (\ref{MTMEqu-0}), using the UL-type of
$RG$-factorization, we obtain%
\[
\left(  I-R_{U}\right)  \left(  I-\Psi_{D}\right)  \left(  I-G_{L}\right)
\mathbf{L}=e\pi-P,
\]
which follows%
\begin{equation}
\left(  I-\Psi_{D}\right)  \left(  I-G_{L}\right)  \mathbf{L}=\left(
I-R_{U}\right)  ^{-1}\left(  e\pi-P\right)  .\label{MTMEqu-1}%
\end{equation}
Let%
\begin{equation}
\mathbf{Y}=\left(  I-G_{L}\right)  \mathbf{L.}\label{MTMEqu-2}%
\end{equation}
Then it follows from (\ref{MTMEqu-1}) that%
\begin{equation}
\left(  I-\Psi_{D}\right)  \mathbf{Y}=\left(  I-R_{U}\right)  ^{-1}\left(
e\pi-P\right)  .\label{MTMEqu-3}%
\end{equation}
We write%
\[
\mathbf{Y}=\left(
\begin{array}
[c]{c}%
\mathbf{Y}_{0}\\
\mathbf{Y}_{1}%
\end{array}
\right)  ,
\]
where $\mathbf{Y}_{0}$ is a matrix with the first $m_{0}$ row vectors of the
matrix $\mathbf{Y}$;%
\[
\left(  I-R_{U}\right)  ^{-1}\left(  e\pi-P\right)  =\left(
\begin{array}
[c]{c}%
\mathbf{H}_{0}\\
\mathbf{H}_{1}%
\end{array}
\right)  ,
\]
where%
\[
\mathbf{H}_{0}=\left(  H_{0,0},H_{0,1},H_{0,2},\ldots\right)  ,
\]%
\[
\mathbf{H}_{1}=\left(
\begin{array}
[c]{cccc}%
H_{1,0} & H_{1,1} & H_{1,2} & \cdots\\
H_{2,0} & H_{2,1} & H_{2,2} & \cdots\\
H_{3,0} & H_{3,1} & H_{3,2} & \cdots\\
\vdots & \vdots & \vdots &
\end{array}
\right)  ,
\]%
\[
H_{i,j}=\left\{
\begin{array}
[c]{ll}%
e\pi_{j}+\sum\limits_{k=1}^{\infty}X_{k}^{\left(  i\right)  }e\pi_{j}%
-X_{j-i}^{\left(  i\right)  }, & j\geq i+1,i\geq0,\\
e\pi_{j}+\sum\limits_{k=1}^{\infty}X_{k}^{\left(  i\right)  }e\pi_{j}-\left[
A_{1}^{\left(  i\right)  }+X_{1}^{\left(  i\right)  }A_{2}^{\left(
i+1\right)  }\right]  , & j=i,i\geq0,\\
e\pi_{j}+\sum\limits_{k=1}^{\infty}X_{k}^{\left(  i\right)  }e\pi_{j}%
-A_{2}^{\left(  i\right)  }, & j=i-1,i\geq1,\\
e\pi_{j}+\sum\limits_{k=1}^{\infty}X_{k}^{\left(  i\right)  }e\pi_{j}, & j\leq
i-2,i\geq2.
\end{array}
\right.
\]
It follows from (\ref{MTMEqu-3}) that%
\begin{equation}
\left(  I-U_{0}\right)  \mathbf{Y}_{0}=\mathbf{H}_{0}\label{MTMEqu-4}%
\end{equation}
and%
\begin{equation}
\left(  I-\Phi_{D}\right)  \mathbf{Y}_{1}=\mathbf{H}_{1}.\label{MTMEqu-5}%
\end{equation}
It follows from (\ref{MTMEqu-5}) that%
\begin{equation}
\mathbf{Y}_{1}=\left(  I-\Phi_{D}\right)  ^{-1}\mathbf{H}_{1},\label{MTMEqu-6}%
\end{equation}
and from (\ref{MTMEqu-4}) that%
\begin{equation}
\mathbf{Y}_{0}=Z\mathbf{H}_{0}+ec_{0},\label{MTMEqu-7}%
\end{equation}
where $c_{0}=\left(  c_{0,0},c_{0,1},c_{0,2},\ldots\right)  $ is an arbitrary
row vector.

Based on (\ref{MTMEqu-6}) and (\ref{MTMEqu-7}), we obtain%
\[
\mathbf{Y}=\left(
\begin{array}
[c]{c}%
\mathbf{Y}_{0}\\
\mathbf{Y}_{1}%
\end{array}
\right)  =\left(
\begin{array}
[c]{c}%
Z\mathbf{H}_{0}+ec_{0}\\
\left(  I-\Phi_{D}\right)  ^{-1}\mathbf{H}_{1}%
\end{array}
\right)  .
\]
It follows from (\ref{MTMEqu-2}) that%
\[
\left(  I-G_{L}\right)  \mathbf{L}=\left(
\begin{array}
[c]{c}%
Z\mathbf{H}_{0}+ec_{0}\\
\left(  I-\Phi_{D}\right)  ^{-1}\mathbf{H}_{1}%
\end{array}
\right)  .
\]
Thus we obtain%
\[
\mathbf{L}=\left(  I-G_{L}\right)  ^{-1}\left(
\begin{array}
[c]{c}%
Z\mathbf{H}_{0}+ec_{0}\\
\left(  I-\Phi_{D}\right)  _{1}^{-1}\mathbf{H}_{1}%
\end{array}
\right)  ,
\]
where%
\[
L_{i,j}=\left\{
\begin{array}
[c]{ll}%
ZH_{0,j}+ec_{0,j}, & i=0,j\geq0,\\
Y_{i}^{\left(  i\right)  }\left(  ZH_{0,j}+ec_{0,j}\right)  +\left(
I-U_{i}\right)  ^{-1}H_{i,j}+\sum\limits_{k=1}^{i-1}Y_{i-k}^{\left(  i\right)
}\left(  I-U_{k}\right)  ^{-1}H_{k,j}, & i\geq1,j\geq0.
\end{array}
\right.
\]

\subsection{The generalized Kemeny's constant}

In this subsection, we generalize the Kemeny's constant (e.g., see Hunter
\cite{Hun:2011}) of the Markov chains from the finite state space to the
infinite state space, where a key for such a generalization is to apply the
censoring technique and the UL-type of $RG$-factorization.

If the QBD $P$ is irreducible and positive recurrent, then the equation
$\left(  I-P\right)  x=0$ exists the unique, up to multiplication by a
positive constant, solution: $x=\gamma e$, where $\gamma$ is positive
constant. Note that such a $\gamma$ is called the Kemeny's constant in the
study of Markov chains with finite states. Here, we extend the Kemeny's
constant to that in a Markov chain with infinite states. That is, the vector
equation (\ref{MTEqu-0}) is shown to have the unique solution%
\begin{equation}
\mathbf{\eta}^{T}=\eta e, \label{MTEqu-1}%
\end{equation}
and $\eta$ is called the generalized Kemeny's constant of the Markov chain
with infinite states, including the QBD process. Using the UL-type of
$RG$-factorization and the censoring technique, we establish a useful relation
between the Kemeny's constant (for the finite states) and the generalized
Kemeny's constant (for the infinite states).

Using the UL-type of $RG$-factorization, it follows from (\ref{MTEqu-0}) that%
\[
\left(  I-R_{U}\right)  \left(  I-\Psi_{D}\right)  \left(  I-G_{L}\right)
\mathbf{\eta}^{T}=0,
\]
which follows%
\begin{equation}
\left(  I-\Psi_{D}\right)  \left(  I-G_{L}\right)  \mathbf{\eta}%
^{T}=0.\label{MTEqu-2}%
\end{equation}
Let%
\begin{equation}
\mathbf{\xi}^{T}=\left(  I-G_{L}\right)  \mathbf{\eta}^{T}.\label{MTEqu-3}%
\end{equation}
Then it follows from (\ref{MTEqu-2}) that%
\begin{equation}
\left(  I-\Psi_{D}\right)  \mathbf{\xi}^{T}=0.\label{MTEqu-4}%
\end{equation}
We write%
\[
\mathbf{\eta}^{T}=\left(
\begin{array}
[c]{c}%
\mathbf{\eta}_{0}\\
\mathbf{\eta}_{1}\\
\mathbf{\eta}_{2}\\
\vdots
\end{array}
\right)  ,\text{ \ \ \ }\mathbf{\xi}^{T}=\left(
\begin{array}
[c]{c}%
\mathbf{\xi}_{0}\\
\mathbf{\xi}_{1}\\
\mathbf{\xi}_{2}\\
\vdots
\end{array}
\right)  .
\]
It follows from (\ref{MTEqu-4}) that%
\begin{equation}
\left(  I-U_{0}\right)  \mathbf{\xi}_{0}=0\label{MTEqu-5}%
\end{equation}
and for $k\geq1$%
\begin{equation}
\left(  I-U_{k}\right)  \mathbf{\xi}_{k}=0.\label{MTEqu-6}%
\end{equation}
Since the matrix $I-U_{k}$ is invertible for $k\geq1$, it is clear from
(\ref{MTEqu-6}) that $\mathbf{\xi}_{k}=0$ for $k\geq1$. Note that the Markov
chain $U_{0}$ is irreducible and positive recurrent, thus we obtain%
\[
\mathbf{\xi}_{0}=\delta e,
\]
where $\delta$ is the Kemeny's constant of the censoring Markov chain $U_{0}$.
Thus we obtain%
\begin{equation}
\mathbf{\xi}^{T}=\left(
\begin{array}
[c]{c}%
\delta e\\
0\\
0\\
\vdots
\end{array}
\right)  .\label{MTEqu-7}%
\end{equation}
It follows from (\ref{MTEqu-3}) and (\ref{MTEqu-7}) that%
\[
\left(  I-G_{L}\right)  \mathbf{\eta}^{T}=\left(
\begin{array}
[c]{c}%
\delta e\\
0\\
0\\
\vdots
\end{array}
\right)  .
\]
This gives%
\begin{align*}
\mathbf{\eta}^{T} &  =\left(  I-G_{L}\right)  ^{-1}\left(
\begin{array}
[c]{c}%
\delta e\\
0\\
0\\
\vdots
\end{array}
\right)  =\left(
\begin{array}
[c]{ccccc}%
I &  &  &  & \\
Y_{1}^{\left(  1\right)  } & I &  &  & \\
Y_{2}^{\left(  2\right)  } & Y_{1}^{\left(  2\right)  } & I &  & \\
Y_{3}^{\left(  3\right)  } & Y_{2}^{\left(  3\right)  } & Y_{1}^{\left(
3\right)  } & I & \\
\vdots & \vdots & \vdots & \vdots & \ddots
\end{array}
\right)  \left(
\begin{array}
[c]{c}%
\delta e\\
0\\
0\\
\vdots
\end{array}
\right)  \\
&  =\delta\left(
\begin{array}
[c]{c}%
e\\
Y_{1}^{\left(  1\right)  }e\\
Y_{2}^{\left(  2\right)  }e\\
\vdots
\end{array}
\right)  =\delta\left(
\begin{array}
[c]{c}%
e\\
e\\
e\\
\vdots
\end{array}
\right)  =\delta e,
\end{align*}
since $Y_{k}^{\left(  k\right)  }e=e$ for $k\geq1$ according to the fact that
$G_{k}e=e$ for $k\geq1$ if the QBD process is irreducible and positive
recurrent, e.g., see Li \cite{Li:2010}. Note that $\mathbf{\eta}^{T}=\eta e$
and $\mathbf{\eta}^{T}=\delta e$, we have $\eta=\delta$. It follows from
(2.19) in Hunter \cite{Hun:2006} that%
\begin{equation}
\eta=\delta=\text{tr}\left(  Z\right)  =\text{tr}\left(  \left(
I-U_{0}+ev_{0}\right)  ^{-1}\right)  .\label{Constant-1}%
\end{equation}

Now, we further apply the censoring technique to computing the generalized
Kemeny's constant $\eta$ of the QBD process. Let%
\[
U_{0}=\left(
\begin{array}
[c]{cc}%
P_{1,1} & P_{1,2}\\
P_{2,1} & P_{2,2}%
\end{array}
\right)  ,
\]
where the sizes of the two matrices $P_{1,1}$ and $P_{2,2}$ are $2$ and
$m_{0}-2$, respectively. Then we can obtain a new censoring chain%
\[
U_{0}^{\left\{  1,2\right\}  }=P_{1,1}+P_{1,2}\left(  I-P_{2,2}\right)
^{-1}P_{2,1}=\left(
\begin{array}
[c]{cc}%
1-a & a\\
b & 1-b
\end{array}
\right)  ,
\]
where $0<a,b\leq1$, since the censoring Markov chain $U_{0}^{\left\{
1,2\right\}  }$ is irreducible and positive recurrent. Let $\delta_{1,2}$ be
the Kemeny's constant of the Markov chain $U_{0}^{\left\{  1,2\right\}  }$.
Then it follows from (3.1) in Hunter \cite{Hun:2006} that%
\begin{equation}
\eta=\delta_{1,2}=1+\frac{1}{a+b}. \label{Constant-2}%
\end{equation}

It is seen from (\ref{Constant-1}) and (\ref{Constant-2}) that the censoring
technique and the UL-type of $RG$-factorization can be applied to highlight
the Kemeny's constant of the Markov chains from the finite state space to the
infinite state space, and also provide an effective algorithm for computing
the generalized Kemeny's constant in the Markov chains with infinite states.

\section{The Variance of the Mixing Time}

In this section, for the QBD process we apply the UL-type of $RG$%
-factorization as well as the generalized inverse to computing the variances
of the first passage time and of the mixing time. Note that the variances can
be expressed by means of the $R$- and $G$-measures through the computational
procedure given in Subsection 2.4.

\subsection{The variance of the first passage time}

To compute the matrix $\mathbf{M}^{\left(  2\right)  }$ given by the first
passage time, we need to solve the equation (\ref{M2Equ-0}). Using the UL-type
of $RG$-factorization, we obtain%
\[
\left(  I-R_{U}\right)  \left(  I-\Psi_{D}\right)  \left(  I-G_{L}\right)
\mathbf{M}^{\left(  2\right)  }=\left\{  E+P\left\{  2\mathbf{M}-\left[
\text{diag}\left(  \pi\right)  \right]  ^{-1}\left[  I+2\left(  e\pi
\mathbf{M}\right)  _{d}\right]  \right\}  \right\}  ,
\]
which follows%
\begin{equation}
\left(  I-\Psi_{D}\right)  \left(  I-G_{L}\right)  \mathbf{M}^{\left(
2\right)  }=\left(  I-R_{U}\right)  ^{-1}\left\{  E+P\left\{  2\mathbf{M}%
-\left[  \text{diag}\left(  \pi\right)  \right]  ^{-1}\left[  I+2\left(
e\pi\mathbf{M}\right)  _{d}\right]  \right\}  \right\}  .\label{M2Equ-1}%
\end{equation}
Let%
\begin{equation}
\mathbf{X}=\left(  I-G_{L}\right)  \mathbf{M}^{\left(  2\right)
}.\label{M2Equ-2}%
\end{equation}
Then%
\begin{equation}
\left(  I-\Psi_{D}\right)  \mathbf{X}=\left(  I-R_{U}\right)  ^{-1}\left\{
E+P\left\{  2\mathbf{M}-\left[  \text{diag}\left(  \pi\right)  \right]
^{-1}\left[  I+2\left(  e\pi\mathbf{M}\right)  _{d}\right]  \right\}
\right\}  .\label{M2Equ-3}%
\end{equation}
We write%
\[
\mathbf{X}=\left(
\begin{array}
[c]{c}%
\mathbf{X}_{0}\\
\mathbf{X}_{1}%
\end{array}
\right)  ,
\]
where $\mathbf{X}_{0}$ is a matrix with the first $m_{0}$ row vectors of the
matrix $\mathbf{X}$;%
\[
\left(  I-R_{U}\right)  ^{-1}\left\{  E+P\left\{  2\mathbf{M}-\left[
\text{diag}\left(  \pi\right)  \right]  ^{-1}\left[  I+2\left(  e\pi
\mathbf{M}\right)  _{d}\right]  \right\}  \right\}  =\left(
\begin{array}
[c]{c}%
\mathbf{S}_{0}\\
\mathbf{S}_{1}%
\end{array}
\right)  ,
\]
where%
\[
\mathbf{S}_{0}=\left(  S_{0,0},S_{0,1},S_{0,2},\ldots\right)  ,,
\]%
\[
\mathbf{S}_{1}=\left(
\begin{array}
[c]{cccc}%
S_{1,0} & S_{1,1} & S_{1,2} & \cdots\\
S_{2,0} & S_{2,1} & S_{2,2} & \cdots\\
S_{3,0} & S_{3,1} & S_{3,2} & \cdots\\
\vdots & \vdots & \vdots &
\end{array}
\right)  ,
\]%
\[
S_{i,j}=E_{i,j}+\Gamma_{i,j}+\sum\limits_{k=1}^{\infty}X_{k}^{\left(
i\right)  }\left(  E_{k+i,j}+\Gamma_{k+i,j}\right)  ,\text{ }i\geq0,j\geq0,
\]%
\[
\Gamma_{i,j}=\left\{
\begin{array}
[c]{ll}%
2\Lambda_{i,j}-A_{0}^{\left(  i\right)  }\left[  \text{diag}\left(  \pi
_{j}\right)  \right]  ^{-1}\left[  I+2\left(  \sum\limits_{i=0}^{\infty}%
e\pi_{i}M_{i,j}\right)  _{d}\right]  , & j=i+1,i\geq0,\\
2\Lambda_{i,j}-A_{1}^{\left(  i\right)  }\left[  \text{diag}\left(  \pi
_{j}\right)  \right]  ^{-1}\left[  I+2\left(  \sum\limits_{i=0}^{\infty}%
e\pi_{i}M_{i,j}\right)  _{d}\right]  , & j=i,i\geq0,\\
2\Lambda_{i,j}-A_{2}^{\left(  i\right)  }\left[  \text{diag}\left(  \pi
_{j}\right)  \right]  ^{-1}\left[  I+2\left(  \sum\limits_{i=0}^{\infty}%
e\pi_{i}M_{i,j}\right)  _{d}\right]  , & j=i-1,i\geq1,\\
2\Lambda_{i,j}, & \text{otherwise},
\end{array}
\right.
\]%
\[
\Lambda_{i,j}=\left\{
\begin{array}
[c]{ll}%
A_{0}^{\left(  0\right)  }M_{1,j}+A_{1}^{\left(  0\right)  }M_{0,j}, &
i=0,j\geq0,\\
A_{0}^{\left(  i\right)  }M_{i+1,j}+A_{1}^{\left(  i\right)  }M_{i,j}%
+A_{2}^{\left(  i\right)  }M_{i-1,j}, & i\geq1,j\geq0.
\end{array}
\right.
\]
It follows from (\ref{M2Equ-3}) that%
\begin{equation}
\left(  I-U_{0}\right)  \mathbf{X}_{0}=\mathbf{S}_{0}\label{SolvM2-1}%
\end{equation}
and%
\begin{equation}
\left(  I-\Phi_{D}\right)  \mathbf{X}_{1}=\mathbf{S}_{1}.\label{SolvM2-2}%
\end{equation}
It follows from (\ref{SolvM2-2}) that%
\begin{equation}
\mathbf{X}_{1}=\left(  I-\Phi_{D}\right)  ^{-1}\mathbf{S}_{1},\label{SolvM2-3}%
\end{equation}
and from (\ref{SolvM2-1}) that%
\begin{equation}
\mathbf{X}_{0}=Z\mathbf{S}_{0}+ec_{0},\label{SolvM2-4}%
\end{equation}
where $c_{0}$ is an arbitrary row vector.

Based on (\ref{SolvM2-3}) and (\ref{SolvM2-4}), we obtain%
\[
\mathbf{X}=\left(
\begin{array}
[c]{c}%
\mathbf{X}_{0}\\
\mathbf{X}_{1}%
\end{array}
\right)  =\left(
\begin{array}
[c]{c}%
Z\mathbf{S}_{0}+ec_{0}\\
\left(  I-\Phi_{D}\right)  ^{-1}\mathbf{S}_{1}%
\end{array}
\right)  .
\]
It follows from (\ref{M2Equ-2}) that%
\[
\left(  I-G_{L}\right)  \mathbf{M}^{\left(  2\right)  }\mathbf{=}\left(
\begin{array}
[c]{c}%
Z\mathbf{S}_{0}+ec_{0}\\
\left(  I-\Phi_{D}\right)  ^{-1}\mathbf{S}_{1}%
\end{array}
\right)  .
\]
Thus we obtain%
\[
\mathbf{M}^{\left(  2\right)  }\mathbf{=}\left(  I-G_{L}\right)  ^{-1}\left(
\begin{array}
[c]{c}%
Z\mathbf{S}_{0}+ec_{0}\\
\left(  I-\Phi_{D}\right)  ^{-1}\mathbf{S}_{1}%
\end{array}
\right)  ,
\]
where%
\[
M_{i,j}^{\left(  2\right)  }=\left\{
\begin{array}
[c]{ll}%
ZS_{0,j}+ec_{0,j}, & i=0,j\geq0,\\
Y_{i}^{\left(  i\right)  }\left(  ZS_{0,j}+ec_{0,j}\right)  +\left(
I-U_{i}\right)  ^{-1}S_{i,j}+\sum\limits_{k=1}^{i-1}Y_{i-k}^{\left(  i\right)
}\left(  I-U_{k}\right)  ^{-1}S_{k,j}, & i\geq1,j\geq0.
\end{array}
\right.
\]

\subsection{The variance of the mixing time}

From (\ref{L2}), we have $\mathbf{L}^{\left(  2\right)  }=\mathbf{M}^{\left(
2\right)  }$diag$\left(  \pi\right)  $. Once the matrix $\mathbf{M}^{\left(
2\right)  }$ is given in Subsection 5.1, it is clear that the matrix
$\mathbf{L}^{\left(  2\right)  }$ is obtained by $\mathbf{M}^{\left(
2\right)  }$diag$\left(  \pi\right)  $.

On the other hand, the matrix $\mathbf{L}^{\left(  2\right)  }$ can be solved
by the matrix equation (\ref{L2Equ-0}) of itself. It may be necessary to
simply provide the outline of solution to Equation (\ref{L2Equ-0}).

To compute the matrix $\mathbf{L}^{\left(  2\right)  }$ given by the mixing
time, we need to solve the equation (\ref{L2Equ-0}). Using the UL-type of
$RG$-factorization, we obtain%
\begin{equation}
\left(  I-\Psi_{D}\right)  \left(  I-G_{L}\right)  \mathbf{L}^{\left(
2\right)  }=\left(  I-R_{U}\right)  ^{-1}\left\{  e\pi+P\left\{
2\mathbf{M}-\left[  \text{diag}\left(  \pi\right)  \right]  ^{-1}\left[
I+2\left(  e\pi\mathbf{M}\right)  _{d}\right]  \right\}  \text{diag}\left(
\pi\right)  \right\}  .\label{L2Equ-1}%
\end{equation}
Let%
\begin{equation}
\mathbf{Y}=\left(  I-G_{L}\right)  \mathbf{L}^{\left(  2\right)
}.\label{L2Equ-2}%
\end{equation}
Then it follows from (\ref{L2Equ-1}) that%
\begin{equation}
\left(  I-\Psi_{D}\right)  \mathbf{Y}=\left(  I-R_{U}\right)  ^{-1}\left\{
e\pi+P\left\{  2\mathbf{M}-\left[  \text{diag}\left(  \pi\right)  \right]
^{-1}\left[  I+2\left(  e\pi\mathbf{M}\right)  _{d}\right]  \right\}
\text{diag}\left(  \pi\right)  \right\}  .\label{L2Equ-3}%
\end{equation}
We write%
\[
\mathbf{Y}=\left(
\begin{array}
[c]{c}%
\mathbf{Y}_{0}\\
\mathbf{Y}_{1}%
\end{array}
\right)  ,
\]
where $\mathbf{Y}_{0}$ is a matrix with the first $m_{0}$ row vectors of the
matrix $\mathbf{Y}$;%
\[
\left(  I-R_{U}\right)  ^{-1}\left\{  e\pi+P\left\{  2\mathbf{M}-\left[
\text{diag}\left(  \pi\right)  \right]  ^{-1}\left[  I+2\left(  e\pi
\mathbf{M}\right)  _{d}\right]  \right\}  \text{diag}\left(  \pi\right)
\right\}  =\left(
\begin{array}
[c]{c}%
\mathbf{T}_{0}\\
\mathbf{T}_{1}%
\end{array}
\right)  ,
\]
where%
\[
\mathbf{T}_{0}=\left(  T_{0,0},T_{0,1},T_{0,2},\ldots\right)  ,,
\]%
\[
\mathbf{T}_{1}=\left(
\begin{array}
[c]{cccc}%
T_{1,0} & T_{1,1} & T_{1,2} & \cdots\\
T_{2,0} & T_{2,1} & T_{2,2} & \cdots\\
T_{3,0} & T_{3,1} & T_{3,2} & \cdots\\
\vdots & \vdots & \vdots &
\end{array}
\right)  ,
\]
and for $i\geq0,j\geq0$
\[
T_{i,j}=e\pi_{j}+\Gamma_{i,j}\text{diag}\left(  \pi_{j}\right)  +\sum
\limits_{k=1}^{\infty}X_{k}^{\left(  i\right)  }\left[  e\pi_{j}%
+\Gamma_{k+i,j}\text{diag}\left(  \pi_{j}\right)  \right]  ,
\]
It follows from (\ref{L2Equ-3}) that%
\begin{equation}
\left(  I-U_{0}\right)  \mathbf{Y}_{0}=\mathbf{T}_{0}\label{L2Equ-4}%
\end{equation}
and%
\begin{equation}
\left(  I-\Phi_{D}\right)  \mathbf{Y}_{1}=\mathbf{T}_{1}.\label{L2Equ-5}%
\end{equation}
It follows from (\ref{L2Equ-5}) that%
\begin{equation}
\mathbf{Y}_{1}=\left(  I-\Phi_{D}\right)  ^{-1}\mathbf{T}_{1},\label{L2Equ-6}%
\end{equation}
and from (\ref{L2Equ-4}) that%
\begin{equation}
\mathbf{Y}_{0}=Z\mathbf{T}_{0}+ec_{0},\label{L2Equ-7}%
\end{equation}
where $c_{0}=\left(  c_{0,0},c_{0,1},c_{0,2},\ldots\right)  $ is an arbitrary
row vector.

Based on (\ref{L2Equ-6}) and (\ref{L2Equ-7}), we obtain%
\[
\mathbf{Y}=\left(
\begin{array}
[c]{c}%
\mathbf{Y}_{0}\\
\mathbf{Y}_{1}%
\end{array}
\right)  =\left(
\begin{array}
[c]{c}%
Z\mathbf{T}_{0}+ec_{0}\\
\left(  I-\Phi_{D}\right)  ^{-1}\mathbf{T}_{1}%
\end{array}
\right)  .
\]
It follows from (\ref{L2Equ-2}) that%
\[
\left(  I-G_{L}\right)  \mathbf{L}^{\left(  2\right)  }=\left(
\begin{array}
[c]{c}%
Z\mathbf{T}_{0}+ec_{0}\\
\left(  I-\Phi_{D}\right)  ^{-1}\mathbf{T}_{1}%
\end{array}
\right)  .
\]
Thus we obtain%
\[
\mathbf{L}^{\left(  2\right)  }=\left(  I-G_{L}\right)  ^{-1}\left(
\begin{array}
[c]{c}%
Z\mathbf{T}_{0}+ec_{0}\\
\left(  I-\Phi_{D}\right)  ^{-1}\mathbf{T}_{1}%
\end{array}
\right)  ,
\]
where%
\[
L_{i,j}^{\left(  2\right)  }=\left\{
\begin{array}
[c]{ll}%
ZT_{0,j}+ec_{0,j}, & i=0,j\geq0,\\
Y_{i}^{\left(  i\right)  }\left(  ZT_{0,j}+ec_{0,j}\right)  +\left(
I-U_{i}\right)  ^{-1}T_{i,j}+\sum\limits_{k=1}^{i-1}Y_{i-k}^{\left(  i\right)
}\left(  I-U_{k}\right)  ^{-1}T_{k,j}, & i\geq1,j\geq0.
\end{array}
\right.
\]

In the rest of this section, we compute $\left(  \mathbf{\eta}^{\left(
2\right)  }\right)  ^{T}$. Note that%
\[
\left(  \mathbf{\eta}^{\left(  2\right)  }\right)  ^{T}=\mathbf{L}^{\left(
2\right)  }e=\mathbf{M}^{\left(  2\right)  }\pi^{T},
\]
hence it is clear that the vector $\left(  \mathbf{\eta}^{\left(  2\right)
}\right)  ^{T}$ is obtained by $\mathbf{M}^{\left(  2\right)  }\pi^{T}$ once
the matrix $\mathbf{M}^{\left(  2\right)  }$ is given in Subsection 5.1.

On the other hand, the vector $\left(  \mathbf{\eta}^{\left(  2\right)
}\right)  ^{T}$ can be solved by the matrix equation (\ref{eta2Equ-0}) of
itself. It may be necessary to simply provide the outline of solution to
Equation (\ref{eta2Equ-0}).

Now, we solve the equation (\ref{eta2Equ-0}). Using the UL-type of
$RG$-factorization, we obtain%
\begin{equation}
\left(  I-\Psi_{D}\right)  \left(  I-G_{L}\right)  \left(  \mathbf{\eta
}^{\left(  2\right)  }\right)  ^{T}=\left(  I-R_{U}\right)  ^{-1}\left\{
e+P\left\{  2\mathbf{M}-\left[  \text{diag}\left(  \pi\right)  \right]
^{-1}\left[  I+2\left(  e\pi\mathbf{M}\right)  _{d}\right]  \right\}  \pi
^{T}\right\}  .\label{etaEqu-1}%
\end{equation}
Let%
\begin{equation}
\Re=\left(  I-G_{L}\right)  \left(  \mathbf{\eta}^{\left(  2\right)  }\right)
^{T}.\label{etaEqu-2}%
\end{equation}
Then it follows from (\ref{etaEqu-1}) that%
\begin{equation}
\left(  I-\Psi_{D}\right)  \Re=\left(  I-R_{U}\right)  ^{-1}\left\{
e+P\left\{  2\mathbf{M}-\left[  \text{diag}\left(  \pi\right)  \right]
^{-1}\left[  I+2\left(  e\pi\mathbf{M}\right)  _{d}\right]  \right\}  \pi
^{T}\right\}  .\label{etaEqu-3}%
\end{equation}
We write%
\[
\Re=\left(
\begin{array}
[c]{c}%
\Re_{0}\\
\Re_{1}%
\end{array}
\right)  ,
\]
where $\Re_{0}$ is a vector with the first $m_{0}$ entries of the vector $\Re
$;%
\[
\left(  I-R_{U}\right)  ^{-1}\left\{  e+P\left\{  2\mathbf{M}-\left[
\text{diag}\left(  \pi\right)  \right]  ^{-1}\left[  I+2\left(  e\pi
\mathbf{M}\right)  _{d}\right]  \right\}  \pi^{T}\right\}  =\left(
\begin{array}
[c]{c}%
\mathbf{W}_{0}\\
\mathbf{W}_{1}%
\end{array}
\right)
\]%
\[
\mathbf{W}_{1}=\left(
\begin{array}
[c]{c}%
\mathbf{W}_{1,1}\\
\mathbf{W}_{2,1}\\
\mathbf{W}_{3,1}\\
\vdots
\end{array}
\right)  ,\text{ }\mathbf{W}_{0}\overset{\text{def}}{=}\mathbf{W}_{0,1},
\]
and for $i\geq0$
\[
\mathbf{W}_{i,1}=e+\sum_{j=0}^{\infty}\Gamma_{i,j}\pi_{j}+\sum\limits_{k=1}%
^{\infty}X_{k}^{\left(  i\right)  }\left[  e+\sum_{j=0}^{\infty}\Gamma
_{k+i,j}\pi_{j}\right]  .
\]
It follows from (\ref{etaEqu-3}) that%
\begin{equation}
\left(  I-U_{0}\right)  \Re_{0}=\mathbf{W}_{0}\label{etaEqu-4}%
\end{equation}
and%
\begin{equation}
\left(  I-\Phi_{D}\right)  \Re_{1}=\mathbf{W}_{1}.\label{etaEqu-5}%
\end{equation}
It follows from (\ref{etaEqu-5}) that%
\begin{equation}
\Re_{1}=\left(  I-\Phi_{D}\right)  ^{-1}\mathbf{W}_{1},\label{etaEqu-6}%
\end{equation}
and from (\ref{etaEqu-4}) that%
\begin{equation}
\Re_{0}=Z\mathbf{W}_{0}+ec_{0},\label{etaEqu-7}%
\end{equation}
where $c_{0}$ is an arbitrary constant.

Based on (\ref{etaEqu-6}) and (\ref{etaEqu-7}), we obtain%
\[
\Re=\left(
\begin{array}
[c]{c}%
\Re_{0}\\
\Re_{1}%
\end{array}
\right)  =\left(
\begin{array}
[c]{c}%
Z\mathbf{W}_{0}+ec_{0}\\
\left(  I-\Phi_{D}\right)  ^{-1}\mathbf{W}_{1}%
\end{array}
\right)  .
\]
It follows from (\ref{etaEqu-2}) that%
\[
\left(  I-G_{L}\right)  \left(  \mathbf{\eta}^{\left(  2\right)  }\right)
^{T}\mathbf{=}\left(
\begin{array}
[c]{c}%
Z\mathbf{W}_{0}+ec_{0}\\
\left(  I-\Phi_{D}\right)  ^{-1}\mathbf{W}_{1}%
\end{array}
\right)  .
\]
Thus we obtain%
\[
\left(  \mathbf{\eta}^{\left(  2\right)  }\right)  ^{T}\mathbf{=}\left(
I-G_{L}\right)  ^{-1}\left(
\begin{array}
[c]{c}%
Z\mathbf{W}_{0}+ec_{0}\\
\left(  I-\Phi_{D}\right)  _{1}^{-1}\mathbf{W}%
\end{array}
\right)  ,
\]
where%
\[
\left(  \eta_{i}^{\left(  2\right)  }\right)  ^{T}=\left\{
\begin{array}
[c]{ll}%
ZW_{0,1}+ec_{0}, & i=0,\\
Y_{i}^{\left(  i\right)  }\left(  ZW_{0,1}+ec_{0}\right)  +\left(
I-U_{i}\right)  ^{-1}W_{i,1}+\sum\limits_{k=1}^{i-1}Y_{i-k}^{\left(  i\right)
}\left(  I-U_{k}\right)  ^{-1}W_{k,1}, & i\geq1.
\end{array}
\right.
\]

Finally, we obtain%
\[
\left(  \mathbf{V}^{\left(  2\right)  }\right)  ^{T}=\left(  \mathbf{\eta
}^{\left(  2\right)  }\right)  ^{T}-\eta^{2}e.
\]
and%
\begin{align*}
v^{\left(  2\right)  }  &  =\pi\left(  \mathbf{V}^{\left(  2\right)  }\right)
^{T}=\pi\left(  \mathbf{\eta}^{\left(  2\right)  }\right)  ^{T}-\eta^{2}\\
&  =-\eta^{2}+\pi_{0}\left(  ZW_{0,1}+ec_{0}\right)  +\sum_{i=1}^{\infty}%
\pi_{i}\left[  Y_{i}^{\left(  i\right)  }\left(  ZW_{0,1}+ec_{0}\right)
\right. \\
&  +\left.  \left(  I-U_{i}\right)  ^{-1}W_{i,1}+\sum\limits_{k=1}%
^{i-1}Y_{i-k}^{\left(  i\right)  }\left(  I-U_{k}\right)  ^{-1}W_{k,1}\right]
.
\end{align*}

\section{Concluding Remarks}

In this paper, we develop some matrix Poisson's equations satisfied by the
mean and variance of the mixing time in an irreducible positive-recurrent
discrete-time Markov chain with infinitely-many levels, and provide a
computational framework for the solution to the matrix Poisson's equations by
means of the UL-type of $RG$-factorization as well as the generalized
inverses. In an important special case: the level-dependent QBD Processes, we
provide a detailed computation for the mean and variance of the mixing time
through the matrix-analytic method.

The results of this paper can be applied to performance computation of
stochastic models such as a discrete-time $MAP/PH/c$ queue and a discrete-time
$MAP/PH/1$ retrial queue by means of the mixing time. Our future work in this
direction contains several different research lines:

(1) Provide algorithms for computing the solution to the matrix Poisson's equations,

(2) analyzing performance measures of practical stochastic models, and

(3) extend the method of this paper, which is based on the UL-type of
$RG$-factorization as well as the generalized inverses, to study more general
Markov models including, Markov chains of $GI/G/1$ type, level-dependent
Markov chains of $M/G/1$ type and of $GI/M/1$ type, continuous-time
block-structured Markov chains, and Markov renewal processes, e.g., see Li
\cite{Li:2010} and Hunter \cite{Hun:1982}.

\vskip        0.5cm

\section*{Acknowledgements}

The first author acknowledges that this research is partly supported by the
National Natural Science Foundation of China (No. 71271187) and the Hebei
Natural Science Foundation of China (No. A2012203125)..


\begin{thebibliography}{9}                                                                                                %

\bibitem {Ald:1983}D. Aldous. \textit{Random walk on finite groups and rapidly
mixing Markov chains}. Lecture Notes in Mathematics, Vol. 986, pages 243--297,
Springer, New York, 1983.

\bibitem {Ald:2002}D. Aldous and J. Fill. \textit{Reversible Markov chains and
random walks on graphs}. 2002. http://www.stat.berkeley.edu/users/aldous/RWG/book.html.

\bibitem {Ald:1997}D. Aldous, L. Lov\'{a}sz and P. Winkler. Mixing times for
uniformly ergodic Markov chains. \textit{Stochastic Processes and their
Applications}, 71(2): 165--185, November 1997.

\bibitem {Asm:1994}S. Asmussen and M. Bladt. Poisson's equation for queues
driven by a Markovian marked point process. \textit{Queueing Systems}, 17:
235--274, 1994.

\bibitem {Bla:1996}M. Bladt. The variance constant for the actual waiting time
of the PH/PH/1 queue. \textit{The Annals of Applied Probability}, 6(3):
766--777, August 1996.

\bibitem {Bri:1995}L. W. Bright and P. G. Taylor. Calculating the equilibrium
distribution in level-dependent quasi-birth-and-death processes.
\textit{Stochastic Models}, 11(2): 497--526, March 1995.

\bibitem {Bri:1996}L. W. Bright and P. G. Taylor. Equilibrium distribution in
level-dependent quasibirth-and-death processes. In \textit{Matrix Analytic
Methods in Stochastic Models}, pages 359--375, Marcel Dekker, New York, 1996.

\bibitem {Cao:2007}X. R. Cao. \textit{Stochastic Learning and Optimization: A
Sensitivity-Based Approach}. Springer-Verlag: New York, 2007.

\bibitem {Cao:1997}X. R. Cao and H. F. Chen. Perturbation realization,
potentials, and sensitivity analysis of Markov processes. \textit{IEEE
Transaction on Automatic Control}, 42 (10): 1382--1393, October 1997.

\bibitem {Cao:1996}X. R. Cao, X. M. Yuan and L. Qiu. A single sample
path-based performance sensitivity formula for Markov chains. \textit{IEEE
Transaction on Automatic Control}, 41(12): 1814--1817, December 1996.

\bibitem {Cat:2010}M. Catral, S. J. Kirkland, M. Neumann and N. S.
Sze.\textbf{ }The Kemeny constant for finite homogeneous ergodic Markov
chains. \textit{Journal of Scientific Computing}, 45(1--3): 151--166, October 2010.

\bibitem {da:2002}A. da Silva Soares and G. Latouche. The group inverse of
finite homogeneous QBD processes. \textit{Stochastic Models}, 18(1):159--171,
January 2002.

\bibitem {Den:2013}S. Dendievel, G. Latouche and Y. Liu. Poisson's equation
for discrete-time quasi-birth-and-death processes. \textit{Performance
Evaluation}, \textbf{Available in Press}, Online 12 June 2013.

\bibitem {Gly:1994}P. W. Glynn. Poisson's equation for the recurrent M/G/1
queue. \textit{Advances in Applied Probability}, 26: 1044--1062, 1994.

\bibitem {Gly:1996}P. W. Glynn and S. P. Meyn. A Liapounov bound for solutions
of the Poisson equation. \textit{The Annals of Probability}, 24(2): 916--931, 1996.

\bibitem {Hey:1995}D. P. Heyman and D. P. O'Leary. What is fundamental for
Markov chains? First passage times, Fundamental matrices, and group
generalised Inverses. In \textit{Proceedings of the 2nd International Workshop
on the Numerical Solution of Markov Chains}, pages 151--159, Kluwer Academic
Publishers, 1995.

\bibitem {Hun:1982}J. J. Hunter. Generalized inverses and their application to
applied probability problems. \textit{Linear Algebra and its Applications},
45: 157--198, June 1982.

\bibitem {Hun:1988}J. J. Hunter. Characterizations of generalized inverses
associated with Markovian kernels. \textit{Linear Algebra and its
Applications}, 102: 121--142, April 1988.

\bibitem {Hun:2006}J. J. Hunter. Mixing times with applications to perturbed
Markov chains. \textit{Linear Algebra and its Applications}, 417(1): 108--123,
August 2006.

\bibitem {Hun:2008}J. J. Hunter. Variances of First Passage Times in a Markov
chain with applications to Mixing Times. \textit{Linear Algebra and its
Applications}, 429(5--6): 1135--1162, September 2008.

\bibitem {Hun:2009}J. J. Hunter. Coupling and mixing times in a Markov chain.
\textit{Linear Algebra and its Applications}, 430(10): 2607--2621, May 2009.

\bibitem {Hun:2011}J. J. Hunter. The role of Kemeny's constant in properties
of Markov chains. In \textit{Proceedings of MSMPRF 2011, Markov and
semi-Markov Processes and Related Fields}, 2011. Available in \textbf{arXiv
preprint arXiv:1208.4716}.

\bibitem {Hun:2013}J. J. Hunter. The distribution of mixing times in Markov
chains. \textit{Asia-Pacific Journal of Operational Research}, 30(1): 1--29,
February 2013.

\bibitem {Kem:1981}J. G. Kemeny. Generalization of a fundamental matrix.
\textit{Linear Algebra and its Applications}, 38: 193--206, June 1981.

\bibitem {Kem:1960}J. G. Kemeny and J. L. Snell. \textit{Finite Markov
Chains}. Van Nostrand, New York, 1960.

\bibitem {Lat:1999}G. Latouche and V. Ramaswami. \textit{Introduction to
Matrix Analytic Methods in Stochastic Modeling}. SIAM Publishing,
Philadelphia,, 1999.

\bibitem {Lev:2009}D. Levin, Y. Peres and E. Wilmer. \textit{Markov Chains and
Mixing Times}. Americian Mathematical Society, 2009.

\bibitem {Li:2010}Q. L. Li (2010). \textit{Constructive Computation in
Stochastic Models with Applications: The RG-Factorizations}. Springer and
Tsinghua Press.

\bibitem {Li:2013}Q. L. Li (2013). Tail probabilities in queueing processes.
To appear in \textit{Asian-Pacific Journal of Operational Research}, Special
Issue of Retrial Queues. Available in \textbf{arXiv preprint arXiv:1209.5604}.

\bibitem {LiC:2004}Q. L. Li and J. Cao. Two types of $RG$-factorizations of
quasi-birth-and-death processes and their applications to stochastic integral
functionals. \textit{Stochastic Models}, 20(3), 299--340, September 2004.

\bibitem {Li:2004}Q. L. Li and L. M. Liu. An algorithmic approach on
sensitivity analysis of perturbed QBD processes. \textit{Queueing Systems},
48(3--4): 365--397, November 2004.

\bibitem {Li:2002}Q. L. Li, and Y. Q. Zhao. A constructive method for finding
$\beta$-invariant measures for transition matrices of M/G/1 type. In
\textit{Matrix Analytic Methods: Theory and Applications}, pages 237--263,
World Scientific: New Jersey, 2002.

\bibitem {Li:2003}Q. L. Li, and Y. Q. Zhao. $\beta$-invariant measures for
transition matrices of GI/M/1 type. \textit{Stochastic Models}, 19(2):
201--233, March 2003.

\bibitem {Li:2009}N. Li, L. Zheng and Q. L. Li. Performance analysis of
two-loop closed production systems. \textit{Computers \& Operations Research},
36(1): 119--134, January 2009.

\bibitem {Liu:2012}J. Liu, S. Yang, A. Wu and S. J. Hu. Multi-state throughput
analysis of a two-stage manufacturing system with parallel unreliable machines
and a finite buffer. \textit{European Journal of Operational Research},
219(2): 296--304, June 2012.

\bibitem {Lov:1998}L. Lov\'{a}sz and P. Winkler. Mixing times. In
\textit{Microsurveys in Discrete Probability}, pages 185--133, American
Mathematical Society, 1998.

\bibitem {Mak:2002}A. M. Makowski and A. Shwartz. The Poisson equation for
countable Markov chains: Probabilistic methods and interpretations. In
\textit{Handbook of Markov Decision Processes}, pages 269--303, Kluwe Academic
Publishers, Norwell, 2002.

\bibitem {Mey:1975}C. D. Meyer, Jr. The role of the group generalized inverse
in the theory of finite Markov chains. \textit{SIAM Review}, 17(3): 443--464, 1975.

\bibitem {Mey:1993}S. P. Meyn and R. L. Tweedie. \textit{Markov Chains and
Stochastic Stability}. Springer, London, 1993.

\bibitem {Mon:2006}R. Montenegro and P. Tetali. Mathematical aspects of mixing
times in Markov chains. \textit{Foundations and Trends in Theoretical Computer
Science}, 1(3): 237--354, 2006.

\bibitem {Neu:1981}M. F. Neuts. \textit{Matrix-Geometric Solutions in
Stochastic Models-An Algorithmic Approach}. The Johns Hopkins University
Press: Baltimore, 1981.

\bibitem {Neu:1989}M. F. Neuts. \textit{Structured Stochastic Matrices of
M/G/1 Type and Their Applications}. Marcel Decker Inc., New York, 1989.

\bibitem {Num:1991}E. Nummelin. On the Poisson equation in the potential
theory of a single kernel. \textit{Mathematica Scandinavica}, 68(1): 59--82, 1991.

\bibitem {Ram:1996}V. Ramaswami. A tutorial overview of matrix analytic
method: with some extensions \& new results. In \textit{Matrix Analytic
Methods in Stochastic Models}, pages 261--296, Marcel Dekker: New York, 1996.

\bibitem {Rui:2009}J. E. Ruiz-Castro, G. Fern\'{a}ndez-Villodre and R.
P\'{e}rez-Oc\'{o}n. A multi-component general discrete system subject to
different types of failures with loss of units. Discrete Event Dynamic
Systems, 19(1): 31--65, March 2009.

\bibitem {Vig:2013}V. Vigon (2013). LU-factorization versus Wiener-Hopf
factorization for Markov chains. To appear in \textit{Acta Applicandae
Mathematicae}, \textbf{Available in Press}, Online February 2013.

\bibitem {Wan:2010}Y. Wang, C. Lin and Q. L. Li. Performance analysis of email
systems under three types of attacks. \textit{Performance Evaluation}, 67(6):
485--499, June 2010.

\bibitem {Wan:2007}Y. Wang, C. Lin, Q. L. Li and Y. G. Fang. A queueing
analysis for the denial of service (DoS) attacks in computer networks.
\textit{Computer Networks}, 51(12): 3564--3573, August 2007.
\end{thebibliography}
\end{document}